\documentclass[reqno,oneside,12pt]{amsart}

\usepackage[T1]{fontenc}
\usepackage{times,mathptm}
\usepackage{amssymb,epsfig,verbatim,xypic}

\theoremstyle{plain}

\newtheorem{thm}{Theorem}[section]

\newtheorem{cor}[thm]{Corollary}
\newtheorem{pro}[thm]{Proposition}
\newtheorem{lem}[thm]{Lemma}

\theoremstyle{definition}
\newtheorem{que}{Question}[section]
\newtheorem{eg}[thm]{Example}
\newtheorem{rem}[thm]{Remark}

\newenvironment{thm-A}
{\noindent{\bf Theorem A.--}\it}{\\}

\newenvironment{thm-M}
{\noindent{\bf Main Theorem.}\it}{\\}

\newenvironment{thm-AA}
{\noindent{\bf Theorem A'.}\it}{\\}

\newenvironment{thm-B}
{\noindent{\bf Theorem B.--}\it}{\\}

\newenvironment{thm-C}
{\noindent{\bf Theorem C.--}\it}{\\}

\newenvironment{thm-BP}
{\noindent{\bf Bell-Poonen Theorem.--}\it}{\\}

\newenvironment{thm-BB}
{\noindent{\bf Theorem B'.}\it}

\def\C{\mathbf{C}}

\def\R{\mathbf{R}}
\def\Q{\mathbf{Q}}

\def\Z{\mathbf{Z}}
\def\N{\mathbf{N}}

\def\Id{{\mathsf{Id}}}

\def\ii{{\sf{i}}}

\def\bbP{\mathbb{P}}

 \def\disk{\mathbb{D}}

\def\Aut{{\sf{Aut}}}

\def\PGL{{\sf{PGL}}\,}

\def\GL{{\sf{GL}}\,}

\def\Tr{{\mathrm{tr}}}

\def\SO{{\sf{SO}}\,}

\def\SU{{\sf{SU}}\,}
\def\U{{\sf{U}}\,}

\def\SL{{\sf{SL}}\,}

\def\End{{\sf{End}}\,}

\def\povol{\mathrm{povol}}

\def\htop{h_{\mathrm{top}}}
\def\hpol{{\mathrm{h}}_{\mathrm{pol}}}

\def\cov{\mathrm{Cov}_{\varepsilon}(n)}
\def\covK{\mathrm{Cov}_{\varepsilon}(n; K)}
\def\ccap{\mathrm{Cap}_{\varepsilon}(n)}
\def\sep{\mathrm{Sep}_{\varepsilon}(n)}

\def\T{\mathbb{T}}

\def\vol{\mathrm{Vol}}

 \def\htop{{\mathrm{h_{top}}}}

\def\Aff{{\sf{Aff}}\,}
\def\SU{{\sf{SU}}\,}

\def\Ind{{\text{Ind}}}

\def\Pic{{\text{Pic}}}

\addtocounter{section}{0}             
\numberwithin{equation}{section}       

\begin{document}

\setlength{\baselineskip}{0.54cm}        
%
%
\title[Automorphisms with slow dynamics]
{Automorphisms of compact K\"ahler manifolds  with slow dynamics}
\date{2019/2020}
\author{Serge Cantat and Olga Paris-Romaskevich}
\address{Univ Rennes, CNRS, IRMAR - UMR 6625, F-35000 Rennes, France}
\email{serge.cantat@univ-rennes1.fr, olga@pa-ro.net}
%
%

%
%

%
%

\begin{abstract} 
We study the automorphisms of compact K\"ahler manifolds having \emph{slow dynamics}. By adapting Gromov's classical argument, we give an upper bound on the polynomial entropy 
and study its possible values in dimensions $2$ and $3$.  We prove that every automorphism with sublinear derivative growth is an isometry ; a counter-example is given in the $C^{\infty}$ context, answering negatively a question of Artigue, Carrasco-Olivera and Monteverde on polynomial entropy. Finally, we classify minimal automorphisms in dimension $2$ and prove they exist only on tori. We conjecture that this is true for any dimension.
\end{abstract}

\maketitle

\section{Introduction}

\subsection{Automorphisms}
Let $X$ be a  compact K\"ahler manifold of dimension $k$. 
By definition, a holomorphic diffeomorphism   $f\colon X\to X$ is an {\bf{automorphism}} of $X$; the group $\Aut(X)$ of all automorphisms
is a (finite dimensional) complex Lie group, with possibly infinitely many connected components. Its neutral component will be denoted $\Aut(X)^0$; its Lie algebra is the algebra of holomorphic vector fields on $X$. 
Our goal is to study automorphisms whose dynamical behavior has ``low complexity''. The main topics will be:
\begin{itemize}
\item polynomial entropy;
\item growth rate  of the derivative $\parallel Df^n\parallel$;
\item equicontinuity of $(f^n)$ on large open subsets;
\item automorphisms acting minimally, or with no periodic orbit. 
\end{itemize}
We focus on low dimensional manifolds $X$, and state a few conjectures in higher dimension. 

\subsection{Polynomial entropy}\label{par:intro-pol-ent}
Let $(X,d)$ be a compact metric space and $f: X \rightarrow X$ a continuous map. The \textbf{Bowen metric}, at time $n$ for the map $f$, is the distance defined by
 the formula 
\begin{equation}\label{eq:Bowen_metric}
d_n^f(x,y)=\max_{0 \leq j \leq n-1} d(f^j(x), f^j(y)).
\end{equation} 
The $(n, \varepsilon)${\bf{-covering number}} $\cov$ is the minimal number of balls of radius $\varepsilon$ in the metric $d_n^f$ that cover $X$. The
{\bf{topological entropy}}  $\htop(f)\in \R_+\cup\{+\infty\}$ is the double limit 
\begin{equation}\label{eq:topological_entropy}
\htop(f):= \lim_{\varepsilon \rightarrow 0} \limsup_{n \rightarrow \infty} \frac{1}{n} \log \left(\cov\right).
\end{equation}
It measures the exponential growth rate of the number of orbits that can be distinguished at a given precision $\epsilon$ during a period of observation equal to $n$.
We are interested in the understanding of "simple" maps, with a topological entropy equal to zero. In this setting, we consider the {\bf{polynomial entropy}} \begin{equation}\label{eq:polynomial_entropy}
\hpol(f):= \lim_{\varepsilon \rightarrow 0} \limsup_{n \rightarrow \infty} \frac{1}{\log n} \log \left(\cov\right).
\end{equation}
This quantity must be taken in $\R_+\cup \{+\infty\}$, but it will be finite for most of the systems we shall consider. 
The polynomial entropy has already been studied in several contexts: for integrable Hamiltonian systems by Marco \cite{Marco:2013,  Marco:2018}, for Brouwer homeomorphisms by Hauseux and Le~Roux \cite{Hauseux-LeRoux}, and in various geometric situations by Bernard, Labrousse and Marco \cite{Bernard-Labrousse, Labrousse:Flat, Labrousse:Geodesic,  Labrousse:Circle, Labrousse-Marco}. A similar notion was defined by Katok and Thouvenot, see \cite{Katok-Thouvenot} and \cite{Kanigowski, KKR}. 

Our first result gives an upper bound on the polynomial entropy of any automorphism $f\colon X\to X$ in terms of the action of $f$ on the cohomology $H^*(X;\C)$; a version of this result is also given for birational transformations. We refer to Sections  \ref{sec:upper_bound} and \ref{subs:proof_birational_theorem_1} for precise statements, and in particular to Theorems \ref{thm:hpol_upper_bound_automorphisms} and  \ref{thm:hpol_upper_bound_birationnal}.
Section \ref{sec:pol_ent_small_dim} gives a refined bound when $\dim(X)$ is small (see Theorem~\ref{thm:dimension_2_and_3_upper_bound}). In Section  \ref{sec:automorphisms_of_surfaces} we are interested in finding all possible values of the polynomial entropy for automorphisms of surfaces: we obtain a partial result in Theorem \ref{thm:rest-entropy}, but we still do not know if there are automorphisms of 
compact K\"ahler surfaces with polynomial entropy in $\R_+\setminus\{0,1,2\}$  (see Question \ref{que:main_two_dimenstions}).

\subsection{Growth of derivatives}
In the setting of ${\mathcal{C}}^\infty$ diffeomorphisms $g\colon M\to M$ of compact manifolds, the growth of the derivatives, i.e. the growth of the 
sequence 
\begin{equation}
\parallel Dg^n\parallel =\max_{x\in M}\parallel D(g^n)_x\parallel,
\end{equation}
where the norm is computed with respect to some fixed riemannian metric on $M$,
 can be slow, for instance less than $n^\alpha$ as $n$ goes to $+\infty$ for every  $\alpha >0$ (see Remark~\ref{rem:Borichev}). In Section \ref{sec:small_entropy}, we describe such an example: it is a variation on classical ideas due to Furstenberg (see Theorem~\ref{thm:Borichev}). As an application, 
 we give a counter-example to a question by Artigue, Carrasco-Olivera and Monteverde \cite{Artigue-Carrasco-Monteverde}: we construct a map with zero polynomial entropy which is not equicontinuous.   
In contrast, Section \ref{sec:slow_growth} shows that automorphisms of compact K\"ahler manifolds do not exhibit such behaviors: 
if the growth of $\parallel Df^n\parallel$ is sublinear, then $f$ preserves a K\"ahler metric (Theorem \ref{thm:derivative_bound}).
This section also studies automorphisms with a dense, Zariski open domain of equicontinuity.

\subsection{Minimal actions}
Another kind of low complexity is when ``all orbits look the same''; one (naive) way to phrase it is to assume
that the action is minimal: every orbit is dense for the euclidean topology. In the smooth setting, there are diffeomorphisms of compact manifolds (and homeomorphisms of surfaces) with positive topological entropy acting minimally. We don't know any such example among automorphisms of compact K\"ahler manifolds. In Section \ref{sec:no_finite_orbits}, we obtain a classification of automorphisms of surfaces satisfying one of the following density properties: $f$ has no finite orbit;  all orbits of $f$ are Zarisky dense; all orbits of $f$ are dense in the euclidian topology ($f$ is minimal).  For instance, minimal automorphisms of surfaces exist only on tori. We conjecture that this is true in any dimension, the first non-trivial case being Calabi-Yau manifolds of dimension $3$.

\subsection{Acknowledgement}

Thanks to Benoit Claudon, Christophe Dupont,  Bassam Fayad, S\'ebastien Gou\"ezel, Jean-Pierre Marco, and Fr\'ed\'eric Le Roux  for interesting discussions on this topic. This work is related to our joint paper with Junyi Xie in which we study free actions of non-abelian free groups~\cite{SOJ}: 
we are really grateful to Junyi Xie for sharing his ideas with us. 

\medskip

\begin{center}
{\bf{ Part I.-- Polynomial entropy : upper bound for automorphisms}}
\end{center}

\section{Upper bound on the polynomial entropy}\label{sec:upper_bound}

Let $X$ be a compact K\"ahler manifold  of (complex) dimension $k$. Let $f$ be an automorphism of $X$. 
Its action on the cohomology of 
$X$ provides a linear map $f^*\colon H^*(X;\Z)\to H^*(X;\Z)$ that preserves the Dolbeault cohomology groups $H^{p,q}(X;\C)\subset H^*(X;\C)$;
we denote by $f^*_j$  the action of $f$ on $H^{j,j}(X;\C)$ or $H^{j,j}(X;\R)$. By definition, the {\bf{polynomial growth rate}} $s_j(f)\in \R_+\cup\{+\infty\}$ 
is the  number
\begin{equation}
s_j(f):=\lim_{n\to +\infty} \frac{\log \parallel (f^{n})_j^*\parallel }{\log(n)}.
\end{equation}
We shall see that $s_j(f)$ is a non-negative integer when the topological entropy of $f$ is equal to $0$, and is infinite if 
$\htop(f)>0$. Then, we set  
\begin{equation}\label{eq:sf}
s(f)=\sum_{j=0}^{\dim_\C(X)} s_j(f)=\sum_{j=1}^{\dim_\C(X)-1} s_j(f).
\end{equation}

\begin{thm}
\label{thm:hpol_upper_bound_automorphisms}
Let $X$ be a compact K\"ahler manifold. If  $f$ is an automorphism of $X$ with $\htop(f)=0$, then $\hpol(f)$ 
is finite and is bounded from above by the following integers
\begin{align*}
  \dim_\C(X) + s(f), \quad 
\dim_\C(X)( s_1(f)+1), \quad
 \dim_\C(X) \times b_2(X).
\end{align*}
\end{thm}

This theorem is the main goal of this section.
Its proof follows Gromov's original argument providing an upper bound for the topological entropy of a holomorphic
endomorphism. 

\begin{rem}
In~\cite{LoBianco}, Lo Bianco proved that the sequence $j\mapsto s_j(f)$ is concave. Since $s_0(f)=0$, and $s_{k-j}(f)=s_j(f)$, this implies 
\begin{equation}
0\leq s_j(f)\leq \min\{ j, k-j\} \times s_1(f).
\end{equation} We shall see that $s_j(f)\leq b_{2j}(X)-1$ if $\htop(f)=0$. 
\end{rem}

\subsection{Gromov's upper bound}\label{subs:Gromov's_upper_bound}

\begin{thm}[\cite{Gromov:entropy}, \cite{Yomdin:2,Yomdin:1}]\label{thm:Gromov_Yomdin}
Let $f: X \rightarrow X$ be a holomorphic endomorphism of a  compact K\"ahler manifold  $X$. 
Then $\htop (f) = \log \lambda (f)$, where $\lambda(f)$ is the spectral radius of the action of $f^*$ on the cohomology $H^*(X, \C)$.
\end{thm}
In fact, Yomdin proved the lower bound $\log \lambda(f)\leq \htop(f)$ for ${\mathcal C}^\infty$ maps of compact manifolds, and Gromov
obtained the upper bound $\htop (f) \leq \log \lambda (f)$ for holomorphic transformations of compact K\"ahler manifolds. 
Gromov's proof first relates the topological entropy to the volumes of {\bf{iterated graphs}} $\Gamma(n)$, and then bounds them by a cohomological 
computation; both steps make use of the K\"ahler assumption. Here,  $\Gamma(n)$ is the image of $X$ under the map $x \mapsto (x, f(x), \ldots, f^{n-1}(x))$. That is,
\begin{equation}
\Gamma(n):=\left\{
\boldsymbol{x}=(x_0, x_1, \ldots, x_{n}) \in X^{n+1} \; |\quad  x_j=f(x_{j-1})
\right\},
\end{equation}
with $n\in \Z_+^*$.
The iterated graphs $\Gamma(n)$ are subsets of $X^{n+1}$, and  $X^{n+1}$ is endowed with the  distance
\begin{equation}
d^X_n (\boldsymbol{x}, \boldsymbol{x}'):=\max_{0 \leq j \leq n} d(x_j, x_j') 
\end{equation}
for every pair of points $\boldsymbol{x}$ and $\boldsymbol{x}'$ in $X^{n+1}$ (as in Section~\ref{par:intro-pol-ent}, here $d$ is some fixed distance on $X$).
Let $\varepsilon$ be a positive real number. By definition, the {\bf{ $(n,\varepsilon)$-capacity}} $\ccap$ is the minimal number of balls of radius $\varepsilon$ in the metric $d^X_n$ that cover $\Gamma(n) \subset X^{n+1}$. A set $S$ is $(n,\varepsilon)$-{\bf{separated}} if $d_n^X(x,y)>\varepsilon$ for every pair of elements $x\neq y$ in $S$, and the $(n,\varepsilon)$-{\bf{separation constant}} $\sep$ is the maximal number of elements in such a set. 

\begin{lem}\label{prop:compare}
For all $n\geq 1$ and $\epsilon >0$ we have 
\begin{align*}
 \ccap = \cov \quad {\text{and}} \quad
\mathrm{Sep}_{2\varepsilon}(n) \leq  \ccap \leq \mathrm{Sep}_{\varepsilon}(n).
\end{align*}
\end{lem}

In particular, one can replace $\cov$ by $\ccap$ or by $\sep$ in the definition of the topological and polynomial entropies without changing their values. 


\begin{proof}
For $\boldsymbol{x}$ and  $\boldsymbol{x}'$ in  $\Gamma(n)$, one has $d_n^X(\boldsymbol{x}, \boldsymbol{x}')=d_n^f(x_0,x_0')$, and the first equality follows. The comparison between $\mathrm{Sep}$ and $\mathrm{Cap}$ holds for every metric space, and in particular for $(X,d_n^f)$. Indeed, if $y_1$, $\ldots$, $y_\ell$ 
are $2\varepsilon$-separated, two of them can not be in the same ball of radius $\varepsilon$; this proves $\mathrm{Sep}_{2\varepsilon}(n) \leq  \ccap$. And if 
$\{y_1, \ldots, y_\ell\}$ is a maximal set of $\varepsilon$-separated points, then every point $x$ is at a distance $\leq \varepsilon$ from one of the $y_j$, proving $\ccap \leq \mathrm{Sep}_{\varepsilon}(n)$.
\end{proof}

Denote by $\pi_j\colon X^{n+1}\to X$ the projection on the $j$-th factor, for $j=0, \ldots,n$. Now, fix a K\"ahler metric on $X$, defined by some K\"ahler form $\kappa$, and put the metric on $X^n$ which is
defined by the K\"ahler form $\kappa_n=\sum_j \pi_j^*\kappa$. This metric differs from $d^X_n$ (as $\ell^2$ norm differs from $\ell^\infty$ norm).  
Let $\vol (\Gamma(n))$ be the $2k$-dimensional volume of $\Gamma(n)$ with respect to the metric $\kappa_n$.  

In order to relate $\vol (\Gamma(n))$ to the $(n, \varepsilon)$-capacity, consider the following definition. Let $W$ be a submanifold of $X^n$ of dimension $\dim_\C(W)=d$.
The $\varepsilon${\bf{-density}} $\mathrm{Dens}_{\varepsilon}(W,z)$ at a  point $z \in W$ is the volume of the intersection of $W$ with a $\kappa_n$-ball $B_z(\varepsilon)$ of radius $\varepsilon$  centered at $z$:
\begin{equation}
\mathrm{Dens}_{\varepsilon}(W,z):=\vol_{2d} \left(
W \cap B_z(\varepsilon)
\right).
\end{equation}
The $(\varepsilon,n)$-{\bf{density}}  $\mathrm{Dens}_{\varepsilon}(W)$ is defined as the  infimum
\begin{equation}
\mathrm{Dens}_{\varepsilon}(W):=\inf_{z \in W} \mathrm{Dens}_{\varepsilon}(W,z).
\end{equation}
Set ${\mathrm{Dens}}_{\varepsilon}(n):= {\mathrm{Dens}}_{\varepsilon}(\Gamma(n))$. Then, 
$\sep {\mathrm{Dens}}_{\frac{\varepsilon}{2}}(n) \leq \vol (\Gamma(n))$ for any $\varepsilon>0$ and Lemma \ref{prop:compare} gives
\begin{equation}\label{eq:capacity_bound}
\log \mathrm{Cap}_{2 \varepsilon}(n) \leq \log \vol (\Gamma(n)) - \log \mathrm{Dens}_{\varepsilon}(n).
\end{equation}
Then, Gromov makes two crucial observations. Firstly, complex submanifolds of compact K\"ahler manifolds are locally minimal 
for any K\"ahler metric (Federer's theorem), and this forces a lower bound for the density: 

\begin{thm}\label{thm:density}
Fix a K\"ahler metric $\kappa$ on $X$, and a real number $\varepsilon>0$.
There exists a positive constant $C=C(\varepsilon, \kappa)$ that does not depend on $n$ such that $\mathrm{Dens}_{\varepsilon}(n) \geq C >0$.
\end{thm}

We refer to \cite{Griffiths-Harris}, pages 389 to 392, for a proof of this result (see also \cite{Gromov:entropy} and \cite{Stolzenberg}). 

Thus, if we divide Equation~\eqref{eq:capacity_bound} by $n$ the term $n^{-1}\log \mathrm{Dens}_{\varepsilon}(n)$ becomes negligeable; as a consequence, $\htop(f)\leq \limsup_n n^{-1} \log \vol (\Gamma(n))$.  For polynomial entropy, we divide by $\log(n)$ and obtain
\begin{equation}\label{eq:polynomial_entropy_bound}
\hpol (f) \leq \limsup_{n \rightarrow \infty} \frac{\log {\vol (\Gamma (n))}}{\log n}.
\end{equation}
Note that both  sides of this inequality are infinite when $\htop (f)>0$.

Secondly, Gromov remarks that this volume growth may be estimated by looking at the action of $f$ on the cohomology of $X$. This comes from the definition of $\Gamma(n)$, and
from the fact that the volume of a complex submanifold of a K\"ahler manifold can be computed homologically: $\vol(\Gamma(n))$
is equal to the pairing  of its homology class with the cohomology class of $\kappa_n^{k}$. We reproduce this argument below to obtain Theorem~\ref{thm:hpol_upper_bound_automorphisms}.

\subsection{Proof of Theorem~\ref{thm:hpol_upper_bound_automorphisms}.}\label{subs:proof_thm_1}

Fix an automorphism $f$ of a compact K\"ahler manifold $X$, and assume that $\htop(f)=0$.

\subsubsection{Action on cohomology} 

\begin{lem}\label{lem:unipotent}
Let $M$ be a compact manifold. There is an integer $m>0$, depending only on $\dim(H^*(M;\R))$, 
with the following property.
If $g$ is a ${\mathcal{C}}^\infty$-diffeomorphism of $M$ and $\htop(g)=0$, then all eigenvalues of 
$(g^m)^*\colon H^*(X;\R)\to H^*(X;\R)$ are equal to $1$. In particular, $g^*$ is virtually unipotent.
\end{lem}

\begin{proof}
Since $g^*$ preserves the integral cohomology $H^*(M;\Z)$ its characteristic polynomial $\chi_{g^*}(t)$ is an element of $\Z[t]$ with leading coefficient equal to $1$. Hence, the eigenvalues of $g^*$ are algebraic integers. Yomdin's lower bound $\log \lambda(g) \leq \htop (g)=0$ shows that all roots of $\chi_{g^*}(t)$ have modulus $\leq 1$, and 
by Kronecker lemma they are roots of $1$ (see~\cite{Kronecker}). Since the degree of $\chi_{g^*}(t)$ is equal to $\dim H^*(X;\R)$,
the order of these roots of unity divides some fixed integer $m$ that depends only on $\dim (H^*(X;\R))$. 
This $m$ satisfies the property stated in the lemma. \end{proof}

\begin{rem}
Denote by $b(M)$ the dimension of $H^*(M;\R)$, and by $m_0(M)$ the optimal value of the integer $m$.
In the proof of Lemma~\ref{lem:unipotent}, one could add the argument given by Levitt and Nicolas in~\cite{Levitt-Nicolas} (proof of Proposition~1.1) to 
conclude that the supremum limit of $\log(m_0(M))/\sqrt{b(M)\log(b(M))}$ as $b(M)$ goes to $+\infty$ is $\leq 1$. 
\end{rem}

Lemma~\ref{lem:unipotent} shows that $f^*$ is virtually unipotent. 
Since $\hpol(f^n)=\hpol(f)$ for all $n\neq 0$, we can replace $f$ by  $f^{m_0(X)}$ to assume that $f^*$ is unipotent.
Fix a norm $\parallel \cdot \parallel$ on the cohomology groups of $X$. Since $f^*$ is unipotent, there 
is a basis of $H^{j,j}(X;\C)$ in which $f^*_j$ is a diagonal of Jordan blocks. The number $s_j(f)$ is the polynomial
growth rate of $\parallel (f^n)^*_j\parallel$ and $s_j(f)+1$ is therefore equal to the size of the largest Jordan block of $f_j^*$. As a consequence, 
\begin{equation}
s_j(f)\leq h^{j,j}(X)-1\leq b_{2j}(X)-1.
\end{equation}

 
\subsubsection{Volumes of iterated graphs}
Now, we use one more time that $X$ is a K\"ahler manifold of dimension $k$.
\begin{thm}[Wirtinger]
Let $Y$ be a compact K\"ahler manifold with a fixed K\"ahler form  $\kappa_Y$. 
If $W$ is a complex analytic subset of $Y$ of dimension $d$, its volume with respect to $\kappa_Y$ is 
equal to 
$\vol (W) = \int_W (\kappa_Y)^{d}= [W] \cdot [\kappa_Y]^{d} $.

\end{thm} 

With $W=\Gamma(n)\subset X^{n+1}$, we obtain
\begin{equation}\label{eq:vol-graph}
\vol (\Gamma(n))=\int_{\Gamma(n)} \kappa_n^{k}=\int_{X} \left(
\sum_{j} \pi_j^* \kappa
\right)^k=\int_X \left(
\sum_{j=0}^{n} (f^j)^* \kappa
\right)^k.
\end{equation}

From the previous paragraph, we know that $\parallel (f^n)^*_1[\kappa] \parallel \leq C_1 \parallel [\kappa]\parallel n^{s_1(f)}$ for some uniform constant $C_1>0$. As a consequence, the norm of the class 
$[\sum_{j=0}^n (f^j)^* \kappa ]$ is no more than $C' \parallel \kappa\parallel n^{s_1(f)+1}$, for some $C' >0$, and since the cup product is a continuous multi-linear map,  we get $\vol (\Gamma(n))\leq C'' n^{k(s_1(f)+1)}$ for some $C''>0$. This proves   the second and third upper bounds of Theorem~\ref{thm:hpol_upper_bound_automorphisms}. 

\begin{lem}\label{lem:pol-estimate}
Let $\ell$ be an integer with $0\leq \ell \leq k$. Then 
\[ 
\parallel (f^{n_1})^*[\kappa] \wedge  \ldots \wedge (f^{n_\ell})^*[\kappa] \parallel
\leq 
C \parallel \! [\kappa]\! \parallel^\ell n_\ell^{s_\ell(f)}\prod_{j=1}^{\ell-1} (n_j-n_{j+1})^{s_j(f)} 
\]
for some constant $C>0$ and every sequence of integers $n_1\geq n_2\geq \ldots \geq n_k\geq 0$.
\end{lem}

\begin{proof}
First, for every class $\omega$ in $H^{j,j}(X;\C)$, $\parallel (f^n)_j^*\omega\parallel \leq C_j n^{s_j(f)}\parallel \omega\parallel$ because the norm of the operators $(f^n)_j^*$ on $H^{j,j}(X;\C)$ 
is bounded by $C_j n^{s_j(f)}$ for some positive constant $C_j, j=1, \ldots, k$. Then, to estimate $\parallel (f^{n_1})^*[\kappa]\wedge (f^{n_2})^*[\kappa]\parallel$ we write 
\begin{eqnarray}
\parallel   (f^{n_1})^*[\kappa]\wedge (f^{n_2})^*[\kappa]\parallel & = & \parallel (f^{n_2})^*\left( (f^{n_1-n_2})^*[\kappa]\wedge [\kappa]\right)\parallel\\
& \leq &  C (n_2)^{s_2(f)}  (n_1-n_2)^{s_1(f)}\parallel \kappa \parallel^2.
\end{eqnarray}
Here the constant $C$ is the product of $C_1$, $C_2$, and a constant $D$ such that $\parallel \omega \wedge \kappa\parallel \leq D\parallel \kappa\parallel \parallel \omega\parallel$ for all classes $\omega$ in $H^{1,1}(X;\C)$. This proves the lemma for $\ell=2$, and this argument extends to other values of $\ell\leq k$.
\end{proof}

Now, by recursion on $\ell$, there is a positive constant $B_\ell$ such that 
\begin{equation}
\sum \parallel(f^{n_1})^*[\kappa] \wedge   \ldots \wedge (f^{n_\ell})^*[\kappa]\parallel
\leq 
B_\ell \parallel [\kappa] \parallel^\ell  n^{\ell + s_1(f)+\ldots s_\ell(f)}
\end{equation}
where the sum is over all $\ell$-tuples $(n_i)$ such that ${n\geq n_1\geq n_2\geq \ldots\geq n_\ell\geq 0}$
With $\ell=k$ we obtain  $\vol (\Gamma(n))\leq B' n^{k+ s(f)}$ for some $B'>0$. This gives the first upper bound 
and concludes the proof of  Theorem~\ref{thm:hpol_upper_bound_automorphisms}.

\subsection{A topological lower bound} 

In Yomdin's proof of  $\htop(f)\geq \log\lambda(f)$ (see Theorem \ref{thm:Gromov_Yomdin}), 
error terms appear in the estimates of $\cov$, and Yomdin shows they grow at most as $\exp(\eta n)$ for every $\eta >0$, 
so that these errors become negligeable when compared to the exponential growth rate of $\cov$ when $\htop(f)>0$  (resp. $\lambda(f)>1$). 

But such an error term appears to be too large to
extract a lower bound on $\hpol(f)$ from the growth $n^{s(f)}$ of $(f^*)^n$ on $H^*(X;\C)$ (see Equation~\ref{eq:sf}). 
We do not know whether all ${\mathcal{C}}^\infty$-diffeomorphisms of compact manifolds (resp. automorphisms of compact K\"ahler manifolds) such that $\parallel (f^n)^*\parallel$ is unbounded (equivalently grows at least like $n$) satisfy $\hpol(f)>0$.

On the other hand, one can adapt arguments of Bowen and Manning  to prove the following proposition. 
Let $G$ be a group with a finite symmetric set of generators $S$. If $\varphi$ is an endomorphism of $G$, set 
\begin{equation}
\rho_{S}(\varphi; G):=\limsup_{n \rightarrow \infty} \frac{\log \mathrm{diam} ( (\varphi)^n(S))}{\log n}
\end{equation} 
where the diameter $ \mathrm{diam}(\cdot )$ is computed with respect to the word length $\vert \cdot \vert_S$ defined by $S$ on $G$. 
For instance, if $h$ is an element of $G$ and $\varphi(g)=hgh^{-1}$ is the conjugacy determined by $h$, 
then $\mathrm{diam} (\varphi)^n(S)\leq 2 \vert h^n\vert_S+1$, with equality when $G$ is a non-abelian free group; thus $\rho_{S}(\varphi; G)=1$
when $G$ is a free group and $h\neq 1_G$.  

\begin{pro}
Let $M$ be a compact manifold; let $S$ be a finite symmetric set of generators of the fundamental group $\pi_1(M)$.
If $g\colon M \rightarrow M$ is a homeomorphism then $\hpol (g) \geq \rho_S(g_*; \pi_1(M)) -1$. 
\end{pro}

Here, $g_*$ is the endomorphism of $\pi_1(M)$ determined by $g$: it is only defined up to composition by a conjugacy
(one has to choose a path that connects the base point $x\in M$ to $g(x)$ and, up to homotopy, two distinct paths differ by an
element of $\pi_1(M;x)$). We skip  the proof because it is exactly the same as in Manning's article~\cite{Manning} and 
because we shall not use this proposition.


\section{Generalization for meromorphic transformations}\label{subs:proof_birational_theorem_1}

In this section, we explain how to extend Theorem~\ref{thm:hpol_upper_bound_automorphisms} 
to the case of meromorphic transformations. To do it, we just have
to replace Gromov's argument by a result of Dinh and Sibony, the drawback being the difficulty to estimate the growth 
of the volumes of the iterated graphs for meromorphic transformations. 
\subsection{Growth on cohomology, graphs and entropy}

Let $g$ be a meromorphic  transformation of a compact K\"ahler manifold $X$ of dimension $k$; let $\Ind(g)$ be its indeterminacy locus. 
Denote by $\left(g\right)_j^*$ the linear action of $g$ on the cohomology group $H^{j,j}(X, \C) \subset H^*(X, \C)$
(see~\cite{Dinh-Sibony:ENS, Guedj:panorama} for a definition), and fix a norm on $H^*(X, \C)$. For every $n\geq 0$, denote by $\parallel (g^n)_j^*\parallel$
the norm of the linear transformation $(g^n)_j^*$, define $s_j(g)\in \R_+\cup \{+\infty\}$ by 
\begin{equation}
s_j(g)=\limsup_{n\to +\infty} \frac{\log\parallel (g^n)_j^*\parallel }{\log(n)}
\end{equation}
and set $s(g):=s_1(g)+\ldots +s_k(g)$ (note that $s_0(g)=0$). 
Since $(g^{n+m})_j^*$ does not coincide with $(g^{m})_j^*\circ (g^{n})_j^*$ in general, it is not clear whether this
supremum limit is actually a limit. The $j$-{\bf{th dynamical degree}} is 
\begin{equation}
\lambda_j(g)=\limsup_{n\to +\infty}  \parallel (g^n)_j^*\parallel ^{1/n},
\end{equation}
and this supremum limit is actually a limit (see~\cite{Dinh-Sibony:ENS}). In these definitions of $s_j(g)$ and $\lambda_j(g)$,  we could replace $\parallel (g^n)_j^*\parallel$ 
by $\int_X (g^n)^*(\kappa^j)\wedge \kappa^{k-j}$ for any fixed K\"ahler form on $X$; this would not change the result (see~\cite{Dinh-Sibony:ENS, BacDang}).

\begin{rem} Except in dimension $\leq 2$, we don't know whether $s_j(g)$ is an integer when $\lambda_j(g)=1$; it could a priori be the case that $ \parallel (g^n)_j^*\parallel $ grows like $\exp(\sqrt{n})$ or $n^{\sqrt{3}}$ 
as $n$ goes to $+\infty$. We refer to \cite{Urech} for this type of questions, and to \cite{Dinh-Sibony:ENS, BacDang} for the main properties
of $\parallel (g^n)_j^*\parallel$.\end{rem}

By definition, the {\bf{iterated graph}} $\Gamma(n)$ of $g$ is the closure 
of the set of points 
\begin{equation}
(x, g(x), \ldots, g^n(x))\in X^{n+1}
\end{equation}
such that $x\notin \Ind(g)$, $g(x)\notin \Ind(g)$, $\ldots$, $g^{n-1}(x)\notin \Ind(g)$. 
To define the notions of entropy, we need to take care of the indeterminacy locus $\Ind(g)$. As in~\cite{Dinh-Sibony:ENS, Guedj:panorama}, we simply use $(n,\varepsilon)$-separated sets (as in \S~\ref{subs:Gromov's_upper_bound} and Lemma~\ref{prop:compare}) for orbits avoiding $\Ind(g)$; in this way, we can talk of topological or polynomial entropy.

\subsection{A bound on the polynomial entropy}

The results of Theorem \ref{thm:hpol_upper_bound_automorphisms} can be extended to meromorphic transformations, as follows. 
Note  that  $\lambda_k(g)$ is the topological degree of $g$ when $k=\dim_\C(X)$; 
thus, if the topological degree of $g$ is $\geq 2$, then $s_k(g)=+\infty$. This is the reason why 
$g$ is assumed to be bimeromorphic in the following statement.

\begin{thm}\label{thm:hpol_upper_bound_birationnal}
If  $g$ is a bimeromorphic transformation of a compact  K\"ahler manifold $X$, then
$\hpol(g)
\leq \dim_\C(X)+ s(g)$.
\end{thm}

\begin{proof}[Sketch of Proof] The first steps of Gromov's argument remain valid: if one defines the polynomial growth
of the iterated graphs by 
\begin{equation}
\povol(g):=\limsup_{n\to + \infty} \frac{\log \vol(\Gamma(n))}{\log(n)},
\end{equation}
then  $\hpol(g)\leq \povol(g)$. Our goal is to show that $\povol(g)$ is bounded from above
by $k+s_1(g)+\ldots+s_k(g)$, where $k=\dim_\C(X)$. For this, we simply copy the argument of \cite{Dinh-Sibony:Annals}. 
More precisely, replacing Lemma~2 of \cite{Dinh-Sibony:Annals} by Corollary~1.2 of \cite{Dinh-Sibony:ENS}, the results of \cite{Dinh-Sibony:Annals}
remain valid on compact K\"ahler manifolds. Thus, there is a positive constant $C$, which depends only on the geometry of $X$, 
such that 
\begin{equation}
\parallel f^*T\parallel \leq C \parallel f^*_j\parallel \parallel T\parallel
\end{equation}
for every meromorphic map $f\colon X\to X$ and every closed positive current $T$ of bi-degree $(j,j)$ on $X$.  Here, $\parallel T\parallel$ 
is the mass of $T$, computed with respect to a fixed K\"ahler form $\kappa$: $\parallel T\parallel=\langle T\vert \kappa^{k-j}\rangle$. 
And $f^*T$ is the positive current which is defined on $X\setminus \Ind(f)$ by pull-back; this upper bound on the mass of $f^*T$ 
implies that the extension of $f^*T$ by $0$ on $\Ind(f)$ is a closed positive current, with the same mass: we shall also denote by
$f^*T$ this current. 

Then, as in the proof of Lemma~5 of \cite{Dinh-Sibony:Annals}, or Lemma~\ref{lem:pol-estimate} above, we obtain the following estimate: for every integer $0\leq \ell\leq k$, 
and every decreasing sequence of integers $n_1\geq n_2\geq \ldots \geq n_k$, 
\begin{equation}
\parallel (g^{n_1})^*[\kappa]  \wedge \ldots \wedge (g^{n_\ell})^*[\kappa]\parallel 
\leq C'  \parallel \!\! [\kappa]\!\! \parallel^\ell \parallel (g^{n_\ell})_{\ell}^* \parallel     \prod_{j=1}^{\ell-1} \parallel (g^{n_j-n_{j+1}})_j^* \parallel 
\end{equation}
for some constant $C'>0$. 
By definition, for every $\eta>0$ and for  $m$ larger than some 
integer $m(\eta;j)$ we have $\parallel (g^m)^*_j\parallel\leq m^{s_j(g)+\eta}$. By recursion on $\ell$, we get 
\begin{equation}
\parallel (g^{n_1})^*[\kappa] \wedge  \ldots \wedge (g^{n_\ell})^*[\kappa]\parallel \leq C'' n^{\ell (1+\eta)
+ s_1(g)+\ldots s_\ell(g)} \parallel [\kappa]\parallel
\end{equation}
for some constant $C''$. To deduce the result, take $\ell=k$ and let $\eta$ go to $0$. 
\end{proof}



\medskip 

\begin{center}
{\bf{Part II.-- Dimension 2: polynomial entropy, minimal actions, \\ and Zariski dense orbits}}
\end{center}

%
%
\section{Polynomial entropy in small dimensions}\label{sec:pol_ent_small_dim}
%
%
The goal of this section is to prove the following theorem.
We also study the case of tori.

\begin{thm}
\label{thm:dimension_2_and_3_upper_bound}
Let $X$ be a compact K\"ahler manifold with $\dim_\C(X) \leq 3$. If $f \in \Aut(X)$ satisfies $\htop (f)=0$ then $\hpol(f) \leq \dim_\C(X)^2$.
\end{thm}

\subsection{Curves and surfaces}\label{par:curves-and-surfaces}

If $X$ is a curve, and $f$ is an automorphism of $X$, the action of $f$ on $H^2(X;\C)$ is the identity, and
Theorem~\ref{thm:hpol_upper_bound_automorphisms} provides the upper bound $\hpol(f)\leq \dim_\C(X)=1$. If the genus of $X$ is positive, then $f$ is an isometry (for the euclidean 
or hyperbolic metric), and $\hpol(f)=0$. If the genus of $X$ is $0$, then $f$ is given by a M\"obius transformation,  and either $f$ is   conjugate to an element of ${\sf{PU}}_2(\C)\subset \PGL_2(\C)$ and then $\hpol(f)=0$, or $f$ has a wandering orbit, and then $\hpol(f)=1$ (see Example~\ref{eg:hpol1_homographies} below). In particular, if $X$ is a curve and $\hpol(f)=0$, then $f$ is an isometry for some K\"ahler metric. 

Suppose now that  $f$ is a bimeromorphic transformation of a compact K\"ahler surface $X$. We shall see in Section~\ref{sec:automorphisms_of_surfaces} that either $\lambda_1(f)>1$, or $\lambda_1(f)=1$ and $s_1(f)\in \{ 0, 1, 2\}$. Thus, {\sl{if $f$ is a bimeromorphic transformation 
of a compact K\"ahler surface, either $\lambda_1(f)>1$, or $\hpol(f)\leq 4$}}. Theorem~\ref{thm:dimension_2_and_3_upper_bound} follows
from this statement and Yomdin's theorem when $\dim_\C(X)=2$.

\subsection{Threefolds}
 
Here we use the results of Lo Bianco (Section $6.2$ in \cite{LoBianco:Thesis} as well as \cite{LoBianco:Thesis_published} and Theorem A in \cite{LoBianco}):

\begin{thm}[Lo Bianco]
Let $f: X \rightarrow X$ be an automorphism of a compact K\"ahler manifold $X$ of dimension $3$. Assume that the
action of $f^*$ on $H^2(X, \C)$ is unipotent; then, it has a unique Jordan block of maximal size, this block is  localized in $H^{1,1}(X, \C)$,
and its size $\ell_1$ satisfies $\ell_1\leq 5$. The other Jordan blocks in $H^2(X, \C)$ have size $\leq \frac{\ell_1+1}{2} \leq 3$. 
 \end{thm}

Since, by duality, the action of $f^*$ on $H^{2,2}(X;\C)$ has Jordan blocks of the same size, we obtain $s_1(f)=s_2(f)\leq 4$, and Theorem~\ref{thm:hpol_upper_bound_automorphisms} gives
\begin{equation}
\hpol(f)\leq 3+ 2\times 4=11.
\end{equation} 
We want to improve this inequality to $\hpol(f)\leq 9$, and for that we use one extra ingredient from the 
proof of Lo Bianco's result. Namely, he found a basis $(u_1, \ldots, u_{\ell_1})$ for the maximal Jordan block of $f^*_1$ that satisfies
\begin{enumerate}
\item $f^*u_1=u_1$ and $f^*u_m=u_m+u_{m-1}$ for any $m=2, \ldots, l_1$ (normal form of the Jordan block);
\item $u_1\wedge u_1=u_1\wedge u_2=0$ in $H^{2,2}(X;\C)$. 
\end{enumerate}
From this result, we can now estimate the volume of $\Gamma(n)$: 
\begin{equation}
\vol (\Gamma(n)) =\int_X \left(
\sum_{j=0}^n (f^j)^* \kappa
\right)^{3} \leq 6  \sum_{i\leq j\leq k=0}^{n} \int_X (f^i)^* \kappa \wedge (f^j)^* \kappa \wedge   (f^k)^* \kappa ;
\end{equation}
since the topological degree of $f$ is  $1$, we can set  $j=i+t_1$ and $k=i+t_2$ to obtain 
\begin{equation}
\vol (\Gamma(n)) \leq 6 \sum_{i=0}^{n} \sum_{t_1=0}^{n-i} \sum_{t_2=0}^{n-i} \int_X \kappa \wedge (f^{t_1})^* \kappa \wedge   (f^{t_2})^* \kappa.
\end{equation}

Denote by $\ell_1> \ell_2 > ...$ the sizes of the Jordan blocks of $f^*$ on $H^{1,1}(X;\C)$. Then, represent the K\"ahler form $\kappa$ 
as a linear combination of vectors 
$\kappa=\sum_{i=1}^{\ell_1} \alpha_i u_i+\sum_{m=1}^{M_2} \sum_{j=1}^{s_2} \beta_i^m v_i^m+\ldots$. Here $M_2$ is the number of Jordan blocks of size $\ell_2$ and the vectors $\{v_i^m\}_{i=1}^{s_2}$ form a basis of the corresponding invariant subspaces in $H^{1,1}(X, \C)$. 
Then write out the wedge product $(f^{t_1})^* \kappa \wedge (f^{t_2})^* \kappa$: since $u_1 \wedge u_1= u_1 \wedge u_2 =0$ and $\ell_2\leq 3$, we see 
that the form $(f^{t_1})^* \kappa \wedge (f^{t_2})^* \kappa$ is a polynomial $P(t_1, t_2)$ in $t_1$ and $t_2$ with values in $H^{2,2}(X,\C)$ and of degree at most $6$.
From this, we get an upper bound

\begin{equation}
\vol \Gamma(n) \leq \sum_{i=0}^n \sum_{t_1=0}^{n-i} \sum_{t_2=0}^{n-i}  \int_X\kappa \wedge P_6(t_1,t_2) \leq C' n^9
\end{equation}
for some positive constant $C'$. This shows that {\sl{an automorphism of a compact K\"ahler manifold of dimension $3$ 
has polynomial entropy $\leq 9$ if its topological entropy vanishes}}. 
This concludes  the proof of Theorem~\ref{thm:dimension_2_and_3_upper_bound}.
 
\subsection{Uniform bound on polynomial entropy}\label{par:2questions}

\begin{que}
Let $f$ be an automorphism of a compact K\"ahler manifold $X$ of dimension $k$. 
If  $\htop (f)=0$ does it follow that $\hpol(f) \leq k^2$ ? Is such an upper bound  optimal, in every dimension $k$ ?
\end{que}


 In dimension $2$,  all the examples for which we are able to calculate the entropy have $\hpol(f) \leq 2$: see the discussion in Section \ref{sec:automorphisms_of_surfaces}. The following proposition provides the $k^2$ upper bound when $X$ is a torus. We provide a cohomological proof, and then give a second, more precise statement in Proposition~\ref{pro:hpol-C-torus}. 

\begin{pro}
If $X$ is a complex torus of dimension $k$, and the automorphism $f\colon X \to X$ satisfies $\htop (f)=0$, then $\hpol(f) \leq k^2$.
\end{pro}

\begin{proof}
We fix a K\"ahler form $\kappa$, and we want to bound:
\begin{eqnarray}\label{eq:volumes_calculation_for_tori_proposition}
\vol (\Gamma (n)) & = & \int_X \left(
\sum_{j=0}^n (f^j)^* \kappa
\right)^{k} \\
& \leq & \sum_{i=0}^n
\sum_{t_1=0}^{n-i} \ldots \sum_{t_{k-1}=0}^{n-i}
 [\kappa] \wedge (f^{t_1})^* [\kappa] \wedge \ldots \wedge (f^{t_{k-1}})^*  [\kappa].
\end{eqnarray} 
The automorphism $f$ is acting linearly on the complex torus $X$ by a matrix $A$,  its action on $H^{1,0}(X, \C)$ is given  by the transposed matrix $A^t$, and  on $H^{0,1}(X, \C)$  by the matrix $\bar{A}^t$. If $\htop(f)=0$ then $f^*$ is virtually unipotent, and we can assume that $f^*$ is unipotent. Fix a basis $(u_1, \ldots, u_k)$ of $H^{1,0}(X, \C)$ in which $A^t$ has a canonical Jordan form, and its  biggest Jordan block corresponds to the subspace generated by $(u_1, \ldots, u_{\ell_1})$; in particular $\ell_1\leq k$. Writing  
$[\kappa] = \sum_{m,n} \alpha_{m,n} u_m\wedge \bar{u}_n $, we obtain 
\begin{equation}\label{eq:representation_degree_f_star}
(f^j)^* [\kappa] = \sum_{m,n=0}^k  \alpha_{m,n} p_{m-1}(j) \bar{p}_{n-1} (j) u_m \wedge \bar{u}_n,
\end{equation}
where the $\alpha_{m,n}$ are complex numbers and the $p_{\delta}(j)$ are polynomial functions of degree $\delta$ in the variable $j$. 
 The maximal degree in the right-hand side of \eqref{eq:representation_degree_f_star} is $2(\ell_1-1)$. Since $u_j \wedge u_j=0$ for all $1\leq j\leq k$,  the sum in Equation~\eqref{eq:volumes_calculation_for_tori_proposition}  is bounded by
\begin{eqnarray}
 \vol (\Gamma (n)) & \leq C & 
 \sum_{i=0}^n
\sum_{t_1=0}^{n-i} \ldots \sum_{t_{k-1}=0}^{n-i} t_1^{2(\ell_1-1)} t_2^{2(\ell_1-2)} \ldots t_{k-1}^{2 \cdot 1} \\
& = & C n^{1+\ell_1-1+ \frac{2 \cdot \ell_1(\ell_1-1)}{2}} = Cn^{\ell_1^2} \leq C n^{k^2}
\end{eqnarray}
for some $C>0$.
\end{proof}

Let $F$ be the element of $\SL_k(\Z)$ given by a Jordan block of size $k$, which means that $F(u_1)=u_1$ and $F(u_m)=u_m+u_{m-1}$ for every $2\leq m\leq k$ in the canonical basis $(u_j)$ of $\Z^k$. This transformation induces a diffeomorphism of the torus $\R^k/\Z^k$ (resp. of the torus $(\C/\Lambda)^k$ for every elliptic curve $\C/\Lambda$).

\begin{lem}\label{lem:hpol-real-torus}
The polynomial entropy of the diffeomorphism $F\colon \R^k/\Z^k\to \R^k/\Z^k$ (resp. of $(\C/\Lambda)^k$) is equal to $k(k-1)/2$ (resp. to $k(k-1)$).
\end{lem}

\begin{proof}[Proof of Lemma~\ref{lem:hpol-real-torus}]
First note that in the canonical basis $(u_j)$ we have
\begin{equation}
F^n=\begin{pmatrix}
1&Q_1(n)&Q_2(n)&\ldots&Q_{k-1}(n)\\
0&1&Q_1(n)&\ldots&Q_{k-2}(n)\\
0&0&1&\ldots&Q_{k-3}(n)\\
& & &\ddots& & \\
0&0&0&\ldots&1
\end{pmatrix},
\end{equation}
where each $Q_j(n)$ is a polynomial function  in the variable $n$ such that $Q_j(n) \approx \frac{n^j}{j!}$ up to lower degree terms.  For simplicity, we set $X:=(\R/\Z)^k$ and choose the $\ell_{\infty}$-metric on $X$. Consider the following set of points $S_n \subset X$: 
\[
S_n:=\left\{
\varepsilon
\left(
i_1, \frac{i_2}{Q_1(n)}, \ldots, \frac{i_k}{Q_{k-1}(n)}
\right) \in X \left\vert\right. i_j \in \Z, 0 \leq i_j \leq \left\lfloor Q_{j-1}(n) \varepsilon^{-1} \right\rfloor
\right\}
\]
with $Q_0(n)=1$ for all $n$ by convention.
Then, 
\begin{equation}
|S_n|=\prod_{j=1}^k \left(
\left\lfloor Q_{j-1}(n) \varepsilon^{-1} \right\rfloor+1
\right) \approx \varepsilon^{-k}\left(\prod_{j=1}^{k-1} j!\right)^{-1} n^{1+2+ \ldots+ (k-1)},
\end{equation}
where the last equivalence holds true up to  terms of lower order in $n$.
The Bowen balls (for $d_n^F$) of radius $\varepsilon$ centered at the points of $S_n$ cover $X$ and, at the same time, all the points in the set $S$  belong to different Bowen balls of radius $\epsilon/2$. Thus
$\cov\simeq |S_n|$ and from the definition of $\hpol$ (see Equation \eqref{eq:polynomial_entropy}) we get $\hpol F = 1+2+ \ldots+ (k-1) = k(k-1)/2$.
\end{proof}

\begin{pro}\label{pro:hpol-C-torus}
If $X$ is a complex torus of dimension $k$, and  $f \in \mathrm{Aut} (X)$ satisfies $\htop f=0$, then $\hpol(f)\leq k(k-1)$.\end{pro}

\begin{proof}
Write $X=\C^k/\Lambda$ for some co-compact lattice $\Lambda\subset \C^k$. There is a matrix $A\in \GL_k(\C)$ 
and a vector $B\in \C^k$ such that $f(z)=A(z) + B\mod \Lambda$. The Bowen distance $d_n^f$ does not depend on 
$B$, so we assume $B=0$ for simplicity. Since $\htop(f)=0$, we can replace $f$ by a positive iterate to assume 
that the action of $f$ on $H^1(X;\R)$ is unipotent (see  Lemma~\ref{lem:unipotent}).

Consider $\C^k$ as a real vector space $V_\R$ of dimension $2k$; fixing a basis of $\Lambda$, we  
identify  it to the lattice $\Z^{2k}\subset V_\R\simeq \R^{2k}$ and denote by $V_\Q\simeq \Q^{2k}$ the 
rational subspace $\Lambda \otimes_\Z \Q$. Since $A$ is a unipotent endomorphism of $V_\Q$, there is basis
of $V_\Q$ in which the matrix of $A$ is a diagonal of Jordan blocks. Since the endomorphism is induced by 
a $\C$-linear transformation, the blocks come in pairs of the same sizes, so that the list of sizes can be
written $k_1\geq  k_2 \geq k_3 \ldots$ with $k_{2i+1}=k_{2i+2}$ for every $i\geq 0$. Now, the proof of Lemma~\ref{lem:hpol-real-torus} and
the additivity  $\hpol(g\times h)=\hpol(g)+\hpol(h)$ give
\begin{equation}
\hpol(f)=\sum_{i\geq 0} k_{2i+1} (k_{2i+1}-1).
\end{equation}
Since $\sum_j k_j=2k$ and $a(a-1)+b(b-1)\leq (a+b)(a+b-1)$ for  all pairs of positive integers, we obtain  $\hpol(f)=k(k-1)$. \end{proof}

\section{Automorphisms of surfaces: classification and lower bounds}\label{sec:automorphisms_of_surfaces}

Consider a homeomorphism $f$ of a compact metric space $(X,d)$. If $K$ is a non-empty compact subset
of $X$, compute the covering number $\covK$ of $K$ by balls of radius $\leq \varepsilon$ in the metric $d_n^f$ (see~Equation~\eqref{eq:Bowen_metric}); 
then, define the polynomial entropy for orbits starting in $K$ by the formula
\begin{equation}
\hpol(f; K) =\lim_{\varepsilon\to 0} \limsup \frac{\log (\covK)}{\log(n)}.
\end{equation}
When $K$ is equal to $X$, we recover the polynomial entropy $\hpol(f)$.
Now, if $U$ is an open subset of $X$, we define the polynomial entropy $\hpol(f;U)$ of $f$ in $U$ to 
be the supremum of the polynomial entropies $\hpol(f;K)$ over all (non-empty) compact subsets $K$ of $U$.
Note that $U$ (resp. $K$) is not assumed to be invariant under the action of $f$. 
This entropy is increasing: if $U$ is contained in $V$ then $\hpol(f;U)\leq \hpol(;V)$; in particular, 
$\hpol(f;U)\leq \hpol(f;X)$. 

Our goal in this section is to prove  the following result. 

\begin{thm}\label{thm:rest-entropy}
Let $f$ be an automorphism of a compact K\"ahler surface $X$ such that $\htop(f)=0$. 
There is 
\begin{itemize}
\item a compact  K\"ahler surface  $X_0$,
\item a regular, bimeromorphic map $\eta\colon X\to X_0$, 
\item  an automorphism $f_0$ of $X_0$ such that $\eta \circ f = f_0\circ \eta$, 
\item and a dense, $f_0$-invariant, Zariski open subset $U$ of $X$,
\end{itemize}
 such that $\hpol(f_0;U)\in\{0, 1, 2\}$. Moreover,  $\hpol(f)=0$ if and only if $f$ preserves a K\"ahler metric on $X$. 
\end{thm}

The first reason to focus on this restricted entropy $\hpol(f_0; U)$ is because blowing-up fixed points 
may change the polynomial entropy (see \S~\ref{par:change_under_blow-ups}). There is also a second reason: 
we were not able to compute $\hpol(f)$ for all automorphisms, for instance for most parabolic automorphisms (see~\S~\ref{sec:surfaces}). 

\begin{que}\label{que:main_two_dimenstions}
Let $f\colon X \to X$ be an automorphism a compact K\"ahler surface such that $\htop(f)=0$. 
Is $\hpol(f)$ an element of $\{ 0, 1, 2\}$ ? \end{que} 

\begin{eg}[Entropy versus restricted entropy]\label{eg:hpol1_homographies}
Let $g$ be a homeomorphism of a compact metric space $M$, with at least one wandering point $x\in M$. Then, $\hpol(f)\geq 1$ (\footnote{Fix a wandering point $x$ and a real number $\epsilon >0$ such that the orbit of the ball $B(x;\epsilon)$ form a family of pairwise disjoint subsets. Set $E(n)=\{ f^{-m}(x)\; ; \; 0\leq m\leq n-1\}$. Then, $E(n)$ is $(n,\epsilon)$-separated: for every pair of distinct points in $E(n)$, say $y=f^{i}(x)$ and $z=f^{j}(x)$ with $i<j$, 
there is a time $m\leq n$, namely $m=j$, such that the distance between $f^m(y)$ and $f^m(z)$ is at least $\epsilon$. Since $E(n)$ has $n$ elements, 
one gets $\hpol(f)\geq 1$.}).
Now, consider a linear projective transformation $A$ of $\bbP^1(\C)$ with a north-south dynamics; up to conjugacy, $A[y_0:y_1]=[ay_0:y_1]$
for some complex number $a$ of modulus $\vert a \vert >1$. Then
\begin{itemize}
\item  $\hpol(A)=1$, 
\item $\hpol(A;U)=0$ if $U=\bbP^1(\C)\setminus\{[0:1], [1:0]\}$;
\item $\hpol(A;V)\in \{0,1\}$ for every open subset $V$ of $\bbP^1(\C)$. 
\end{itemize}
If we start with the parabolic homography $B[y_0:y_1]=[y_0+y_1:y_1]$, then $\hpol(B)=1$ and $\hpol(B;U)=0$ is $U$ is an open subset
that does not contain the unique fixed point of $B$. 
\end{eg} 

\subsection{Entropy and invariant fibrations}\label{par:all-automorphisms}

\subsubsection{Hodge decomposition} 
Let $X$ be a compact K\"ahler manifold. Fix a K\"ahler form $\kappa$ on $X$, and denote by $[\kappa]\in H^2(X;\R)$ its cohomology class. First, recall the Hodge decomposition 
\begin{equation}\label{eq:Hodge_decomposition}
H^n(X, \C)\cong \bigoplus_{p+q=n} H^{p,q} (X, \C)
\end{equation}
where $H^{p,q}(X, \C)$ is the subspace of cohomology classes of type $(p,q)$. This decomposition is invariant under the action of $\mathrm{Aut} (X)$. Moreover, $\kappa$ determines an hermitian form  $Q_\kappa: H^2(X, \C) \times H^2(X, \C) \rightarrow \C$,
\begin{align*}
Q_\kappa([\alpha], [\beta]) = [\alpha]\cup [\bar{\beta}] \cup [\kappa]^{k-2} =\int_X \alpha \wedge \bar{\beta} \wedge \kappa^{k-2}.
\end{align*}
According to the Hodge index theorem, the restriction of $Q_\kappa$ to $H^{1,1}(X,\C)$ has signature 
$(1, h^{1,1}(X)-1)$ where $h^{1,1}(X)=\dim H^{1,1}(X,\C)$. 
Its restriction to   $H^{2,0}(X, \C) \oplus H^{0,2}(X, \C)$ is positive definite.

\subsubsection{Surfaces}\label{sec:surfaces}

Assume that $X$ is a surface, then $k=2$, and we denote $Q_\kappa$ by $Q$ because it does not depend on $\kappa$. The signature of $Q$
on $H^{1,1}(X;\R)$ being equal to $(1, h^{1,1}(X)-1)$, it determines a structure of Minkowski space on $H^{1,1}(X;\R)$. 
The classification of isometries of Minkowski spaces and the geometry of surfaces lead to the following three cases (see \cite{Cantat:Milnor}):

\smallskip

(1) $f^*$ is  an {\bf{elliptic}} isometry of $H^{1,1}(X;\R)$, and then there exists a positive iterate $f^m$ of $f$ such that $f^m\in \Aut(X)^0$. In particular, there 
is a holomorphic vector field $\theta$ on $X$ such that $f^m$ is the flow, at time $1$, obtained by integrating~$\theta$.

\smallskip

(2) $f^*$ is a {\bf{parabolic}} isometry of $H^{1,1}(X;\R)$. In that case, there exists a fibration $\pi\colon X\to B$ onto a Riemann surface $B$ whose 
generic fibers are connected and of genus $1$, which is invariant under the action of $f$: there is an automorphism $f_B$ of $B$ such that 
$\pi\circ f = f_B\circ \pi$. Moreover, the growth of $(f^n)^*$ on $H^{1,1}(X;\R)$ is quadratic: $\parallel (f^n)^*\parallel \approx C n^2$ for some
positive constant~$C$.

\smallskip

(3) $f^*$ is a {\bf{loxodromic}} isometry of $H^{1,1}(X;\R)$, $f^*$ has a unique eigenvalue of modulus $>1$ on $H^*(X;\C)$,  this eigenvalue coincides
with $\lambda(f)$, and it is realized on $H^{1,1}(X;\R)$. According to Theorem~\ref{thm:Gromov_Yomdin}, the topological entropy of $f$ is positive, and equal to $\log(\lambda(f))$.

\smallskip

Moreover, in the loxodromic case, $f\colon X\to X$ has infinitely many saddle periodic points, and these periodic points equidistribute towards
the unique $f$-invariant probability measure of maximal entropy (see~\cite{Cantat:Acta, Cantat:Milnor, Dujardin:Duke}). 

\smallskip

We shall say that the automorphism $f$ is elliptic, parabolic, or loxodromic if its action $f^*$ on $H^{1,1}(X;\R)$ is 
respectively elliptic, parabolic, or loxodromic.

\begin{pro}\label{pro:auto_surfaces_basics}
Let $f$ be an automorphism of a compact K\"ahler surface. If $f$ is loxodromic, then $\hpol(f)=+\infty$.
If $f$ is parabolic, then $\hpol(f)\in [2,4]$, and there is an $f$-invariant effective divisor $D\subset X$ such that
$\hpol(f;X\setminus D)$ is equal to $2$. If $f$ is elliptic, then $\hpol(f)\in [0,2]$.
\end{pro}

\begin{eg}\label{eg:quadric}
 An automorphism
 $g \in \PGL_2(\C)$ of $\bbP^1(\C)$ satisfies $\hpol(g)=1$, except when $g$ is (conjugate to) a rotation in which case $\hpol(g)=0$ (see Example~\ref{eg:hpol1_homographies}).
 Now, consider the group of automorphisms of $\bbP^1 \times \bbP^1$; this group contains $\PGL_2(\C)\times \PGL_2(\C)$ 
 as a subgroup of index $2$. By additivity of polynomial entropy for products, we see that $\{0,1,2\}$ is exactly the set of possible
 polynomial entropies for automorphisms of $\bbP^1 \times \bbP^1$. \end{eg}

\begin{eg}\label{eg:parabolic_torus}
Let $E=\C/\Lambda$ be an elliptic curve. Consider the automorphism of  of $E^2$ defined by $h(x,y)=(x+y,y)\mod(\Lambda^2)$. 
Using real coordinates, $h$ is conjugate, by some real analytic diffeomorphism, to the linear diffeomorphism of $\R^2/\Z^2\times \R^2/\Z^2$
given by $H(x_1,y_1, x_2,y_2)=(x_1+y_1,y_1,x_2+y_2,y_2)$, i.e. to the diagonal diffeomorphism $h_1\times h_2$ where 
$h_i(x_i,y_i)=(x_1+y_1,y_1)$. From Lemma \ref{lem:hpol-real-torus} we know that $\hpol(h_i)=1$; hence, $\hpol(h)=2$. \end{eg}

\begin{proof}
If $f$ is loxodromic, $\htop(f)$ is positive, and then $\hpol(f)$ is infinite.
 If $f$ is elliptic, then some positive iterate of $f$ is in  $\Aut(X)^0$; in that case $s_1(f)=s_2(f)=0$ and Theorem \ref{thm:hpol_upper_bound_automorphisms} gives $\hpol(f) \leq 2$. 

If $f$ is parabolic, $\parallel (f^*)^n \parallel$
is quadratic, hence $\hpol(f)\leq 4$ by Theorem~\ref{thm:hpol_upper_bound_automorphisms}. Consider
the $f$-invariant genus $1$ fibration  $\pi\colon X\to B$, and the action $f_B$ induced by $f$ on the base $B$. 

If $f_B$ has finite order, we replace $f$ by some positive iterate $f^m$ to assume $f_B=\Id_B$; then $f$ acts 
by automorphism on each fiber of $\pi$, and $f$ acts by translation on each regular fiber of $\pi$. Then, we 
define $D$ to be the union of all singular fibers of $\pi$. The open set $U=X\setminus D$ is $f$-invariant, and if
$K$ is a compact subset of $U$, its projection $\pi(K)$ can be covered by a finite number of compact disks $K_i$
such that $\pi$ is equivalent to the trivial fibration $\pi_i\colon K_i\times \R^2/\Z^2\to K_i$ above $K_i$: there is
a real analytic diffeomorphism $\psi_i\colon \pi^{-1}(K_i)\to K_i\times \R^2/\Z^2$ such that (1) $\psi_i$ is an affine map
on each fiber and (2) $ \pi = \pi_i\circ \psi_i$ on $\pi^{-1}(U_i)$.
Conjugating $f\colon \pi^{-1}(U_i)\to \pi^{-1}(U_i)$ by $\psi_i$, we obtain a diffeomorphism $f_i$ of $K_i\times \R^2/\Z^2$
that acts by translations: 
\begin{equation}
f_i(b,z)=(b,z+t(b)), \quad \forall (b,z)\in K_i\times \R^2/\Z^2,
\end{equation}
where $b\in K_i\mapsto t(b)\in \R^2/\Z^2$ is a smooth map. It is shown in~\cite{Cantat:Groups} that $t(b)$ is a real
analytic mapping of rank $2$ : the image $t(K_i)$ contains an open subset of $\R^2/\Z^2$. This implies that $\hpol(f_i)$
is equal to $2$. Indeed, $\hpol(f_i)\leq 2$ because $\parallel Df_i^n\parallel$ grows linearly and the base has dimension $2$; 
and $\hpol(f_i)\geq 2$ because the polynomial entropy of $(x,y)\mapsto (x,x+y)$ on $[a,b]\times \R/\Z$ is equal to $1$
for every interval $[a,b]$ with $a< b$.
This argument shows that $\hpol(f;K)=2$, so that $\hpol(f;U)=2$, as desired. 

Now, assume that $ f_B$ has infinite order. Let us prove that {\sl{$X$ is a (compact complex) torus}}.
From Proposition~3.6 of \cite{Cantat-Favre}, the minimal model $X_0$ of $X$ is a torus. More precisely, 
the Kodaira dimension of $X$ vanishes, $X$ has a unique minimal model $X_0$, given by a (canonical)
birational morphism $\eta\colon X\to X_0$, and this minimal model is a torus $X_0=\C^2/\Lambda$. 
The morphism $\eta$ contracts a divisor $D\subset X$ onto a finite subset 
of $X_0$. The uniqueness of the minimal model implies that $\eta$ is $f$-equivariant: there is an automorphism
$f_0$ of $X_0$ such that $\eta\circ f = f_0\circ \eta$. 
On the other hand, $f$ and $f_0$ are two parabolic automorphisms, preserving a unique fibration of genus $1$; 
so this fibration must be $\eta$-invariant, meaning that $\pi\colon X \to B$ factors as $\pi=\pi_0\circ \eta$ for some 
fibration $\pi_0\colon X_0\to B$. The action of $f$ and $f_0$ on $B$ coincide (they are both given by  $f_B$);  
so, the finite set $\eta(D)$ is a finite $f_0$-invariant subset of $X_0$, and $\pi_0(\eta(D))$ is a finite $f_B$-invariant 
subset of $B$. Since all orbits of $f_B$ are infinite, we deduce that $D$ is empty: this means that $\eta$ is 
an isomorphism, and $X$ is equal to the torus $\C^2/\Lambda$. 

Since $X$ is a torus, $f$ is induced by an affine  automorphism of $\C^2$:  $f(x,y)=A(x,y)+(s,t)$ modulo $\Lambda$,
for some virtually unipotent linear map $A\in \GL_2(\C)$. Then, $\hpol(f)=2$, as shown in Section~\ref{par:2questions}. 
\end{proof}

\subsection{Preliminaries for elliptic automorphisms} 

In order to prove Theorem~\ref{thm:rest-entropy}, we now need to study elliptic automorphisms, and since 
$\hpol(f)=\hpol(f^m)$ for every $m\neq 0$, we may and do assume that $f$ is in $\Aut(X)^0$.
We shall use three basic facts, which we collect in this section. 

\subsubsection{Wandering saddle configuration} 
\begin{figure}
\includegraphics[scale=0.6]{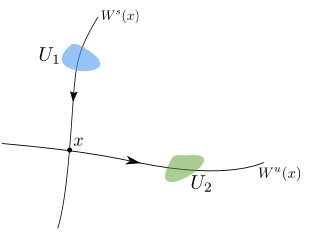}
\caption{Wandering saddle connection of a saddle fixed point $x$.}
\label{fig:wandering_saddle_connection}
\end{figure}

Let $f$ be an automorphism of a compact K\"ahler surface $X$.
We say that $f$ has a
{\bf{wandering saddle configuration}} if $f$ has a saddle fixed point $x$, together with two open subsets $\U_1$ and
$\U_2$ in $X$ such that
\begin{enumerate}
\item[(a)] $\U_1\cap \U_2=\emptyset$;

\item[(b)] $f^n(\U_i)\cap \U_i=\emptyset$ for $i=1,2$ and all $n\neq 0$;

\item[(c)] $\U_1\cap W^s(x)\neq \emptyset$ and $\U_2\cap W^u(x)\neq \emptyset$,
\end{enumerate}
where $W^s(x)$ and $W^u(x)$ denote the stable and unstable manifolds of $x$. See Figure \ref{fig:wandering_saddle_connection} for illustration.
The following lemma is essentially due to Hauseux and Le Roux, \cite{Hauseux-LeRoux}. 

\begin{lem}\label{lem:wsc}
Let $f$ be an automorphism of a complex surface with a wandering saddle configuration. Then $\hpol(f)\geq 2$, and $\hpol(f)=2$ if
$f$ is an element of $\Aut(X)^0$.
\end{lem}

\begin{proof}[Sketch of proof]
Take a point $x_1$ in $\U_1\cap W^s(x)$ and a point $x_2$ in $W^u(x)$.  
Fix $\epsilon>0$ such that, for $i\in \{1,2\}$, the ball of radius $\epsilon$ centered at $x_i$ is contained in $U_i$ and
is at distance $> \epsilon$ from the complement of $\U_i$. If $\epsilon$ is small enough, these balls are wandering.
If $\ell$ is large enough, one can find a point $z_\ell$ in $B(x_1,\epsilon)$ whose orbit $f^n(z_\ell)$ stays in the complement of 
the two balls except for $f^0(z_\ell)\in B(x_1,\epsilon)$ and $f^\ell(z_\ell)\in B(x_2,\epsilon)$. Then, the points $f^{-j}(z_\ell)$
for $j\leq n$ and $\ell\leq n/2$ are $(\epsilon/2,n)$-separated (see \cite{Hauseux-LeRoux} Exemple 2); the size of this
set grows quadratically with $n$, hence $\hpol(f)\geq 2$. The equality follows from $\hpol(f)\leq 2$ when 
$f\in \Aut(X)^0$ (see Proposition~\ref{pro:auto_surfaces_basics}).
\end{proof}

\subsubsection{Compact groups and Kodaira dimension} \label{par:compact_kodaira}
Fix a K\"ahler metric on $X$, given by a K\"ahler form $\kappa_0$. If $f$ is contained in a compact subgroup $K$ of $\Aut(X)$, 
and $\mu$ is a Haar measure on $K$ (with $\mu(K)=1$), the K\"ahler form $\kappa=\int_K g^*\kappa_0\; d\mu(g)$ is $f$-invariant; in particular, 
$f$ preserves a K\"ahler metric if $f$ has finite order. 

\begin{lem}\label{lem:compact_kod_pos}
Let $X$ be a compact K\"ahler manifold. If $\Aut(X)^0$ is not compact, then $X$ is ruled. 
If the Kodaira-dimension $kod(X)$ of $X$ is $\geq 0$, then $\Aut(X)^0$ is a compact Lie group. 
\end{lem}

\begin{proof} This lemma follows from Theorem~4.9 of \cite{Lieberman} and its proof; we only sketch the argument.
From Theorems~3.3 and~3.12 of \cite{Lieberman}, we know that
 $\Aut(X)^0$ is an extension of a compact torus $T_X$ by a linear algebraic group~$L_X$:
\begin{equation}
1\to L_X \to \Aut(X)^0 \to T_X\to 0.
\end{equation}
If $\Aut(X)^0$ is not compact, then $L_X$ is non-trivial and we can fix a one-dimensio\-nal algebraic subgroup $H\subset L_X$. Lieberman
proves that the action of $H$
on $X$ is compactifiable (see \cite{Lieberman}, Chapter 3): this shows that the closures of the orbits of $H$ are rational curves and
that  $X$ is ruled.  Since the Kodaira dimension of a ruled manifold is negative, this concludes the proof of the lemma. 
\end{proof}

We shall say that the iterates of $f$ -- or that $f^\Z$ --  form an {\bf{equicontinuous (or normal) family}} on some open subset $U$ of $X$ if, given 
any subsequence $f^{n_i}$, and any relatively compact open set $V\subset U$, one can extract a further subsequence $f^{n_j}$
that converges uniformly on $V$ towards a continuous map $g\colon V\to X$ (the limit $g$ is automatically holomorphic). 
The equicontinuity is equivalent to a uniform bound on the norm of the differentials $Df^n$ on any relatively compact subset $V$ 
of $U$. If the iterates of $f$ form an equicontinuous family on $X$, then $f^\Z$ is contained in a compact subgroup of $\Aut(X)$ and, by 
Lemma~\ref{lem:compact_kod_pos}, $f$ preserves a K\"ahler form.

Lemma~\ref{lem:compact_kod_pos} shows that $\Aut(X)^0$ preserves a K\"ahler metric when $kod(X)\geq 0$. Thus, in what follows, we assume that 
$kod(X)=-\infty$ and distinguish two cases: 
\begin{itemize}
\item $X$ is an irrational, ruled surface -- this is studied in Section~\ref{par:IRS};
\item $X$ is a rational surface -- this is studied in Sections~\ref{par:plane} and~\ref{par:rational}.
\end{itemize}

\subsubsection{M\"obius twists}  

We say that a function $a\colon {\overline{\disk}} \to \C$ is holomorphic if it extends to a holomorphic function on some 
neighborhood of ${\overline{\disk}}$. We say that $a$ is a holomorphic diffeomorphism onto its image, if this property 
is satisfied by some holomorphic extension.

\begin{lem}\label{lem:skew-P2}
Let $b\colon {\overline{\disk}} \to \C^*$ be a holomorphic function that does not vanish. 
Let $g$ be the transformation of ${\overline{\disk}}\times \bbP^1(\C)$ 
defined by 
$$g(x,[y_0:y_1])=(x,[y_0+b(x) y_1:y_1]).
$$ 
Then, $\hpol(g)=1$. \end{lem}

\begin{proof}
Since $b$ does not vanish, $\varphi(x,[y_0:y_1])=(x,[b(x) y_0 :y_1])$ is a holomorphic diffeomorphism 
of ${\overline{\disk}}\times \bbP^1(\C)$ that conjugates $g$ to $(x,[y_0:y_1])\mapsto (x,[y_0+y_1:y_1])$.
Thus, $\hpol(g)=1$ (see Example~\ref{eg:hpol1_homographies}).
\end{proof}

The following lemma is much more interesting, but will not be used in this article. 
We postpone its proof to the appendix. 

\begin{lem}\label{lem:skew-P1}
Let $a\colon {\overline{\disk}} \to \C^*$ be a holomorphic function such that $\vert a \vert -1$ does not vanish. Let $f$ be the transformation of  $\; {\overline{\disk}}\times \bbP^1(\C)$ 
defined by 
$$f(x,[y_0:y_1])=(x,[a(x)y_0:y_1]).$$
Then, $\hpol(f)=1$. 
\end{lem}

\subsection{Irrational, ruled surfaces}\label{par:IRS}

In this section, $X$ is a ruled surface, but $X$ is not rational. This means that there is a fibration $\pi\colon X\to B$ 
onto a base $B$ of genus $\geq 1$ with generic fiber $\bbP^1$. This fibration is equivariant with respect to $f\in \Aut(X)^0$
and an automorphism $f_B\colon B\to B$.

\subsubsection{}\label{par:change_under_blow-ups} Assume, first, that $f_B$ is periodic. Since the polynomial entropy does not change if one replace $f$ 
by some positive iterate, we may as well suppose that $f_B=\Id_B$. The following lemma does not require the genus of  $B$ 
be positive.

\begin{lem}\label{lem:ruled-irr-surfaces1}
Let $\pi\colon X\to B$ be a ruled surface, and let $f$ be an automorphism of $X$ preserving each fiber of $\pi$. 
Denote by $f_b$ the restriction of $f$ to the fiber $X_b:=\pi^{-1}(b)$, for $b\in B$. Then 
\begin{itemize}
\item either $f$ preserves a K\"ahler metric and $\hpol(f)=0$;
\item or $f_b$ is a loxodromic homography on at least one fiber, and then $\hpol(f)=1$ or $2$;
\item or $f_b$ is a parabolic homography for all fibers except finitely many fibers $F_i$ on which 
$f_b$ is the identity, and then $\hpol(f; V)=1$, where  $V=X\setminus\cup_i F_i$. 
\end{itemize}  
If  $\hpol(f)\geq 1$ there is a Zariski closed, $f$-invariant subset $F\subset X$ such that 
$\hpol(f; X\setminus F)=1$.
\end{lem} 

\begin{proof}[Proof of Lemma~\ref{lem:ruled-irr-surfaces1} when $\pi$ is a submersion] In this first part of the proof, we assume that $\pi$ is a submersion.

\begin{lem}\label{lem:ruled_trick}
Let $\pi\colon X\to B$ be a ruled surface. Assume that $\pi$ is a submersion. 
Let $f$ be an automorphism of $X$ that preserves every fiber of $\pi$. Then, 
\begin{itemize}
\item either the automorphisms $f_b$ of the fibers $\pi^{-1}(b)\simeq \bbP^1(\C)$ are pairwise
conjugate;
\item or, for every $b$, $f_b$ is either parabolic or the identity.
\end{itemize}
\end{lem}

\begin{proof}
For every $b\in B$, the automorphism $f_b$ of the fiber $X_b:=\pi^{-1}(b)$. If we fix an isomorphism $X_b\simeq \bbP^1(\C)$,
$f_b$ becomes a M\"obius transformation of $\bbP^1(\C)$, induced by some matrix $A(b)$ in $\GL_2(\C)$. This matrix is 
not uniquely determined, but the function $b\mapsto  C_f(b)=\Tr(A(b))^2/\det(A(b))$ is well defined because it is invariant 
under scalar multiplication and under conjugacy. 
If $U\subset B$  is a small disk, then on $\pi^{-1}(U)$ the fibration  $(X,\pi)$  is biholomorphically
equivalent to $(U\times \bbP^1(\C),\pi_1)$, where $\pi_1$ denotes the first projection. Using coordinates $(x,[y_0:y_1])$ on $U\times \bbP^1(\C)$, one can write  $f(x,[y_0:y_1])=(x, [A(x)(y_0,y_1)])$ where now $x\mapsto A(x)\in \GL_2(\C)$ is holomorphic. So, the
function $C_f$ is holomorphic. Since $B$ is compact, $C_f$ is constant. 

If the constant $C_b$ is not equal to $4$, the automorphisms $f_b$ are pairwise conjugate: to a rotation of $\bbP^1(\C)$ if 
$C_f\in [0,4[$, to a loxodromic homography otherwise. If $C_b=4$, two cases may occur: $f_b$ can be parabolic, conjugate
to $[y_0:y_1]\mapsto [y_0+y_1:y_1]$, or it can be the identity. \end{proof}

Keeping the same notation as in the proof of Lemma~\ref{lem:ruled_trick}, there are three cases to be distinguished, depending 
on the value of $C_f$.

{\sl{Step 1.-- Rotations}}.  Let us assume that $C_f=(2\cos(2\pi \theta_0))^2$ for some    $\theta_0\neq 0 \mod(1)$.  Then $f$ 
 is locally conjugate to $(x,[y_0:y_1])\mapsto (x,  [\exp(2 i \pi \theta_0) y_0:y_1])$ on $U \times \bbP^1(\C)$. In particular, the iterates of $f$ form an equicontinous family on each open subset $\pi^{-1}(U)$, with $U$ as above. From Section~\ref{par:compact_kodaira}, we know that $f$
preserves a K\"ahler form.

{\sl{Step 2.-- North south dynamics}}. We assume now that $C_f\notin [0,4]$.  Then $f$ is locally conjugate to $(x,[y_0:y_1])\mapsto (x,  [\lambda y_0:y_1])$ 
 is not constant, so that $\hpol(f)\geq 1$ for some $\lambda$ of modulus $\neq 1$. In that case, one gets $\hpol(f)=1$ (and $\hpol(f;V)=0$
 if $V$ is the Zariski open subset obtained by taking of the two sections given by the fixed points of $f$).

{\sl{Step 3.-- Parabolic homographies}}. The last case is $C_f = 4$. If $f$ is not the identity, then 
we define $U_B$ to be the open subset of $B$ obtained by removing the points for which $f_b= \Id_{X_b}$. 
Lemma~\ref{lem:skew-P2} shows that $\hpol(f; \pi^{-1}(U_B))=1$.
\end{proof}

The following argument is not fully necessary to prove Theorem~\ref{thm:rest-entropy}, but it illustrates what may happen when we do a 
blow-up.

\begin{proof}[Proof of Lemma~\ref{lem:ruled-irr-surfaces1}  in the general case] We do not assume that $\pi$ is a submersion anymore; we keep the notation used in the case of a submersion. What may happen is that $\pi$ has a finite number of reducible fibers, made of trees of rational curves. Changing $f$ to a positive 
iterate, we assume that $f$ preserves each irreducible component of each fiber. Contracting some of those
irreducible components, we obtain a birational morphism $\eta\colon X \to X_0$ and a submersion $\pi_0\colon X_0\to B$ 
such that (i) $ \pi =\pi_0\circ \eta$ and (ii) $f$ induces an automorphism $f_0$ of $X_0$ 
with $\eta \circ f = f_0\circ \eta$. There are three cases, according to the value of $C_{f_0}$. 
First, assume that $f_0$ acts by rotations on the fibers of $\pi_0$. If $z\in B$ is a
critical value of $\pi$ and $U$ is a small disk around $z$, then on $\pi_0^{-1}(U)$, $f_0$ is conjugate to 
\begin{equation}\label{eq:example_ruled}
(x,[y_0:y_1])\mapsto (x,  [\lambda y_0:y_1]),
\end{equation}
with $\vert \lambda \vert = 1$. 
To recover $f$, we have to blow-up some fixed points of $f_0$ in the fiber $\pi_0^{-1}(z)$. 
Doing so, we see that, above $U$, the iterates of $f$ also form an equicontinuous family. This implies that $f^\Z$ form an equicontinuous 
family on~$X$, and then that $f$ preserves a K\"ahler form (see Section~\ref{par:compact_kodaira}). 
The second case is when $C_{f_0}\notin [0,4]$, then $f_0$ is also locally conjugate to the map given by Equation~\ref{eq:example_ruled}, but 
with $\vert \lambda\vert \neq 1$. At its two fixed points on $\pi_0^{-1}(z)$, the diferential $Df_0$ is a diagonal map of type $(x,y)\mapsto (x,\alpha y)$
with $\alpha= \lambda^{\pm 1}$. If we blow-up such a fixed point, we create a wandering saddle configuration, with a fixed point at which the 
linear part is $(x,y)\mapsto (x/\alpha, \alpha y)$. Any further blow-up conserves at least one such wandering saddle configuration. So, $\hpol(f)=2$.

The remaining case is when each $(f_0)_{b}$ is either parabolic or the identity. In that case, one can just shrink $U_B$ by taking away the 
critical values of $\pi$; one gets $\hpol(f;\pi^{-1}(U_B))=1$.
\end{proof}

\subsubsection{}\label{par:minimal_ruled_inifinite} Assume that $f_B$ has infinite order. Then, because $g(B)\geq 1$ ($X$ is not rational), we know that $B$ is an elliptic curve
and $f$ is a translation on $B$.
All orbits of $f_B$ and thus of $f$ are infinite, and in particular $\pi$ is a submersion.
In that case, $f$ is the flow, at time $t=1$, of some holomorphic vector field on $X$ that is transverse to the fibration $\pi$; this vector field 
determines a Riccati foliation ${\mathcal{F}}$ on $X$ that is transverse to $\pi$. Let $A$ denote the Zariski closure of $f^\Z$ in the group
$\Aut(X)^0$; then $A$ sits in an extension $1\to K_A\to A \to \Aut(B)^0\to 0$, where $\Aut(B)^0$ can be identified to $B$ (acting on itself
by translation) and $K_A$ is the subgroup of $A$ preserving each fiber of $\pi$. If $\dim(A)=1$, then $K_A$ is finite, $A$ is a compact 
group, $f$ preserves a K\"ahler form and $\hpol(f)=0$. 

Now, assume $\dim(A)\geq 2$, so that the connected component $K_A^0$ of $K_A$ has dimension $\geq 1$.
Consider an element $g$  of $K_A^0\setminus\{ \Id_X\}$, and apply Lemma~\ref{lem:ruled_trick}  to it. 
If $C_g=4$, then $A_g(b)$ is conjugate to 
a parabolic homography $[y_0:y_1]\mapsto [y_0+a_g(b)y_1:y_1]$ for every $b$; since $g$ commutes to $f$, the set of fixed points of $g$ is $f$-invariant,
so it contains no fiber since otherwise $g$ would fix a Zariski dense set of fibers; as a consequence, $a_g(b)$ does not vanish, and the fixed point set
of $g$ coincides with a section $\sigma\colon B\to X$ of $\pi$. 
If $C_g$ is not equal to $4$, then the fixed point set of $g$ intersects every fiber in $2$ points. If the order of $g$ is infinite, which we can assume, 
the derivative of $g_{\vert X_b}$ at each of these fixed points are distinct (and inverse of each other), so that  those two 
fixed points can be distinguished, one from the other, and define two sections of $X$. 

Now, we can apply the arguments of \cite[Proposition 3.1]{Loray-Marin} and of~\cite[Section 3]{Potters:1969} (namely the constructions on pages 251--253) to prove that $X$ and $f$ are given by one of the following three examples. 

\begin{eg}
Up to a finite base change, $X$ is just the product  $B\times \bbP^1(\C)$ and $f(x,y)=(x+\tau, A(y)) $ for some translation $\tau$
and some homography $A$. Then $\hpol(f)\in\{0, 1\}$ and $\hpol(f)=0$ if and only if $f$ preserves a K\"ahler metric (if and only if $A$
is conjugate to a rotation).
\end{eg}

\begin{eg}
There are two sections $\sigma_0$ and $\sigma_\infty$ of the fibration and, if one takes them of, the complement is isomorphic 
to the quotient of $\C^*\times \C^*$ by the action of a cyclic subgroup, generated  by the transformation $\gamma(x,y)\mapsto (\lambda x, \mu y)$, 
with $\vert \lambda\vert<1$ and $B=\C^*/\langle \lambda\rangle$, and some $\mu$ in $\C^*$. The action of $f$ on $X$ lifts to a diagonal transformation $F(x,y)=(\alpha x, \beta y)$ on $\C^*\times \C^*$. Consider the real number $\eta$ such that $\vert \lambda\vert^\eta= \vert \mu\vert$; the function $H(x,y)=
\vert x \vert^\eta/\vert y \vert$ is $\gamma$-invariant: it determines a continuous function $H_X$ on $X$ with zeroes and poles along the two sections. 
Moreover, $H\circ f = \tau H$ with $\tau := \vert \alpha\vert^\eta/ \vert \beta\vert$. If $\tau=1$, then $f$ is contained in a compact subgroup of $\Aut(X)^0$, 
$\hpol(f)=0$ and $f$ preserves a K\"ahler metric. If $\tau\neq 1$, then every orbit of $f$ visits the tube $1\leq H_X\leq \tau$  exactly once, and $\hpol(f)=1$
(because this tube is compact).
\end{eg}

\begin{eg}
There is a unique section, when one removes it the complement is isomorphic to the quotient of $\C^*\times \C$
by $\gamma\colon (x,y)\mapsto (\lambda x, y+1)$ and $f$ lifts to $F(x,y)=(\alpha x, y+\beta)$. Here,  $f$ is contained in a compact
subgroup of $\Aut(X)^0$ if and only if $\beta=\log\vert \alpha\vert / \log\vert \lambda\vert$; otherwise, $\hpol(f)=1$, since 
one can argue as in the previous example, with the  $\gamma$-invariant function  $H(x,y)=y-\log\vert x\vert$.
\end{eg}

Thus, we obtain the following result. 

\begin{lem}\label{lem:ruled-irr-surfaces2}
Let $\pi\colon X\to B$ be a ruled surface over an elliptic curve $B$.
Let $f$ be an automorphism of $X$ such that $\pi\circ f = f_B\circ \pi$ for some automorphism  of $B$ of  infinite order.
Either $f$ is contained in a compact subgroup of $\Aut(X)^0$, or the $\alpha$ and $\omega$-limit set of every orbit
are contained in sections of $\pi$, and in that case, $\hpol(f)=1$.
\end{lem}

\subsubsection{Conclusion for irrational surfaces}\label{par:Irrational_surfaces}
Lemma~\ref{lem:ruled-irr-surfaces1} and~\ref{lem:ruled-irr-surfaces2} prove Theorem~\ref{thm:rest-entropy} when $X$ is a ruled, irrational surface. In the following paragraphs, 
we deal with rational surfaces.

\subsection{Linear projective case}\label{par:plane} Before studying rational surfaces in full generality, we focus on automorphisms
of  $\bbP^2(\C)$. 

\begin{pro}\label{prop:pgl3}
Let $g$ be an element of $\PGL_{3}(\C)=\Aut(\bbP^2(\C))$. Then $\hpol(g) \in \left\{0,1,2\right\}$. More precisely, the following classification holds: 
\begin{enumerate}
\item $\hpol(g)=0$ if and only if $g$ is an isometry;
\item $\hpol(g)=2$ if and only if $g$ has a wandering saddle configuration, if and only if 
 $g$ is conjugate to  
\begin{equation}
\begin{pmatrix}
\lambda&0&0\\
0&\mu&0\\
0&0&1
\end{pmatrix}
\; \; {\text{or}} \; \; 
\begin{pmatrix}
1&1&0\\
0&1&0\\
0&0&\nu
\end{pmatrix}
\end{equation}
with $|\lambda|>1>|\mu|$ and $ |\nu| \neq 1$.
\item $\hpol(g)=1$ if and only if $g$ is conjugate  
\begin{equation}\label{eq:LPC(3)}
\begin{pmatrix}
1&0&0\\
0&\alpha&0\\
0&0&\nu
\end{pmatrix},
\; \; {\text{or}} \; \; 
\begin{pmatrix}
\beta&1&0\\
0&\beta& 0\\
0&0&1
\end{pmatrix},
\; \; {\text{or}} \; \; 
\begin{pmatrix}
1 &1 &0\\
0 & 1 &1\\
0 & 0 & 1
\end{pmatrix}.
\end{equation}
with $ |\nu| \neq 1$ and $|\alpha|=|\beta|=1$.
\end{enumerate}
\end{pro}

In proving this proposition, we also introduce two technics: the blow up of fixed points and the 
symbolic coding of Hauseux and Le Roux. 

\begin{proof} Fix a system of  homogeneous coordinates $[x:y:z]$ of $\bbP^2(\C)$.

The first step is to show that every linear projective transformation $g\in \Aut(\bbP^2(\C))$ is contained in one of the mentionned 
conjugacy classes. This is classical. 

In case (1),  $g$ preserves a K\"ahler metric and its polynomial entropy vanishes. 

In case (2),  $g$ has a wandering saddle configuration: in the diagonal case, on can take the fixed point $q=[0:0:1]$ and the
stable and unstable varieties $\{y=0\}$ and $\{x=0\}$; in the second case one can take $q=[0:1:0]$ and the stable and unstable varieties 
$\{z=0\} $ and $\{ x=0\}$ (assuming $\vert \nu\vert <1$). 

Let us now look at case (3). First, assume that $g$ is (conjugate to) a diagonal transformation with eigenvalues $1$, $\alpha$ and $\nu$ with 
$\vert \alpha\vert=1 < \nu$, i.e. $g(x,y)=(\frac{1}{\nu}x,\frac{\alpha}{\nu} y)$ in affine coordinates $(x,y)$. Blow up the fixed point $[0:0:1]$
to get a new surface $X$ on which $g$  lifts to an automorphism $g_X$: the surface $X$ fibers on $\bbP^1(\C)$: each fiber is 
the strict transform of a line through $[0:0:1]$, the action of $g_X$ preserves this ruling, it acts by a rotation $w\mapsto \alpha w$ 
on the base, and as a loxodromic isometry in the fibers. Let $\Delta$ be a disk in the base of this fibration, centered at one of the fixed 
point ($0$ or $\infty$) of the rotation $w \mapsto \alpha w$. Above $\Delta$, the fibration is holomorphically equivalent to a product 
$\Delta\times \bbP^1(\C)$ on which $g_X$ acts as $(w,[y_0:y_1])\mapsto (\alpha w, [\nu y_0:y_1)$ with $\vert \nu\vert \neq 1$; since 
we can cover the base by two such disks, one sees from Example~\ref{eg:quadric}  that $\hpol(g_X)=1$. Since $g$ 
has a wandering point we obtain $1\leq \hpol(g)\leq \hpol(g_X)=1$, hence $\hpol(g)=1$. 

Assume that, after conjugacy, $g$ is given by the second matrix of the list~\eqref{eq:LPC(3)}. Blow up the fixed point $[1:0:0]$, to
get a new surface $X$, together with a fibration $X\to \bbP^1(\C)$ corresponding to the lines $y=c^{st} z$ of $\bbP^2(\C)$. 
Then $g$ lifts to an automorphism $g_X$ of $X$ preserving the ruling, acting as $w\mapsto \beta w$ on the base, and as $[y_0:y_1]\mapsto [y_0+y_1:y_1]$ 
in the fibers. As in the previous case, we obtain $\hpol(g)=1$.

Assume now that, after conjugacy, $g$ is given by the third matrix of the list~\eqref{eq:LPC(3)}. 
Then $g$ has a unique fixed point $q=[1:0:0]$ and all other points 
are wandering (their $\alpha$ and $\omega$ limit sets coincide with $\{q\}$). This setting has been studied by Hauseux and Le Roux in 
\cite{Hauseux-LeRoux} and we can directly apply their result. Let $X_0$ be the complement of the fixed point $q$. 

Let  $\mathcal{F}=\{ F_1, \ldots F_k\}$ be a finite family of non-empty subsets of $X_0$. Let $F_\infty$ be the complement of $\cup_{F_i\in  \mathcal{F}} F_i$. 
To each orbit $(g^n(x))$, one associates its possible codings, i.e. the sequences of indices $i(n)\in \{1, \ldots, k,\infty\}$ such that $g^n(x)\in F_{i(n)}$ for all $n\in \Z$ (the coding is not unique since the $F_i$ may overlap). Let $Cod(N)$ be the number of codes of length $N$ which are realized by orbits of $g$; the 
polynomial degree growth of $Cod(N)$ is denoted ${\mathrm{Pol_{cod}}}(f;\mathcal{F})$ (this is denoted by $\hpol(f;\mathcal{F})$ in \cite{Hauseux-LeRoux}, \S~2.2). Then, one can define the \textbf{local polynomial entropy of} $g$ \textbf{at a finite subset} $\mathcal{S} \subset X_0$ as the limit of ${\mathrm{Pol_{cod}}}(f;U_n(\mathcal{S}))$ where $U_n(\mathcal{S})$ is any decreasing sequence of open subsets of $X_0$ such that $\cap U_n(\mathcal{S})=\mathcal{S}$ (see \cite{Hauseux-LeRoux}). This number is denoted by $\hpol^{{\mathrm{loc}}}(f;\mathcal{S})$.
One says that  subsets $U_1, \ldots, U_L \subset X_0$ are {\bf{mutually singular}} if for every $M\geq 1$ there exists a point $x \in X$ and times $n_1, \ldots, n_L$ such that
\begin{itemize}
\item[(i)] $g^{n_i}(x) \in U_i$ for every $i$ 
\item[(ii)] $|n_i-n_j|>M$ for every $i \neq j$. 
\end{itemize}
A finite set $\{s_1, \ldots, s_n\} \in  X_0$ is {\bf{singular}} if all small enough neighborhoods $U_1, \ldots, U_L$ of $s_1, \ldots, s_L$ respectively are mutually singular. Any singleton is a singular set.
Hauseux and Le Roux obtain the equality
\begin{equation}
\hpol(f)=\sup \left\{
\hpol^{\mathrm{loc}}(f;\mathcal{S})\; |\; \; \mathcal{S} \; \; \textit{is a singular set}
\right\}.
\end{equation}
Now, coming back to the example of the Jordan block of size $3$, one sees that every singular set of $g$ is reduced to a singleton. Indeed, any singular set 
contains $q$ because $q$ is the unique $\omega$-limit point; and, for any point $p \neq q$, the number $N$ of iterations for which a neighbourhood $U_1$ of $p$ reaches a fixed neighbourhood $U_2$ of $q$ (with  $g^N (U_1) \cap U_2 \neq \emptyset$), is finite. 
This shows that $\hpol(g)=1$.
\end{proof}


\begin{rem}
Note that Proposition \ref{prop:pgl3}, and its proof can be repeated word by word for $\PGL_3 (\R)$. 
\end{rem}

\subsection{Rational surfaces}\label{par:rational} To prove Theorem~\ref{thm:rest-entropy} it remains to study rational surfaces which are not isomorphic to the
projective plane. First, we study minimal, rational surfaces.

\subsubsection{Minimal rational surfaces}\label{par:minimal_rational_surfaces} 
Let $X$ be a minimal rational surface, and assume that $X$ is not the projective plane. 
When $X$ is isomorphic to $\bbP^1\times \bbP^1$, we know from Example~\ref{eg:quadric} that the conclusion of Theorem~\ref{thm:rest-entropy} is satisfied.
Thus, we assume that $X$ is not isomorphic to  $\bbP^1(\C)\times \bbP^1(\C)$ or $\bbP^2(\C)$. Then, there is a unique ruling $\pi \colon X\to \bbP^1(\C)$,
invariant under the action of $\Aut(X)$, with a unique section $C\subset X$, of negative self intersection $C^2=-d$ for some $d\geq 1$ (see~\cite{BPVdV}, \S~V.4).

Denote by $B\simeq \bbP^1(\C)$ the base of the fibration, fix any homogeneous coordinate $[x_0:x_1]$ on $B$, and set $x=x_0/x_1$,  the corresponding 
affine coordinate. Above $x\neq 0$ and $x\neq \infty$, the fibration $\pi$ is equivalent to the trivial fibration $\C\times \bbP^1(\C)\to \C$, with 
the section $C$ corresponding to the point at infinity in $\bbP^1(\C)$, i.e. to $y=\infty$ in the affine coordinate $y=y_0/y_1$ of $\bbP^1(\C)$. 
These two charts are glued together by the map $(x,y)\mapsto (1/x, x^d y)$. Every automorphism $f$ of $X$ can be written in the first 
(resp. second) chart as 
\begin{equation}
f(x,[y_0:y_1])=(f_B(x), [a y_0 + q(x)y_1:y_1])
\end{equation}
for some linear projective transformation $f_B\in \Aut(B)$, some complex number $a\neq 0$, and some polynomial function $q\in \C[x]$ 
of degree at most $d$; going to the second chart replaces $\alpha$ by its inverse, and $q(x)$ by $x^dq(1/x)$.  (More precisely, $q$ should be seen as a homogeneous polynomial of degree $d$ in $(x_0,x_1)$.)

\smallskip

{\bf{A.--}} First, we assume that $f_B(x)=\alpha x$ for some $\alpha\neq 0$ and some appropriate choice of 
 coordinate on $B$. 

\smallskip
 
{\bf{A1.--}} Assume $\vert \alpha\vert \neq 1$.
The fiber $x=0$ is $f$-invariant, and $f$ has one or two fixed points along this fiber.
If $y\mapsto ay+q(0)$ is conjugate to a $y\mapsto y+1$ or to $y\mapsto ay$ with $\vert a \vert \neq 1$, then  $f$ has
a wandering saddle configuration at one of these fixed points and $\hpol(f)=2$ by Lemma~\ref{lem:wsc}. 
So, when $\vert \alpha\vert \neq 1$ we may assume $\vert a \vert =1$; in that case $\alpha^m\neq a$ for all $m\in \Z$.

\smallskip

{\bf{A2.--}} Now, assume that $\alpha^m\neq a$ for all $m\in \{0, 1, \ldots, d\}$, write $q(x)=\sum_k q_k x^k$ and consider the polynomial function 
\begin{equation}
p(x)=\sum_k q_k (a-\alpha^k)^{-1}x^k.
\end{equation}
Then, $h(x,y)=(x,y+p(x))$ defines an automorphism of $X$ (because the degree of $p$ is $\leq d$) that  conjugates $f$ to the automorphism $g(x,y)=(\alpha x, ay)$. We obtain $\hpol(f;U)=1$ if $\vert \alpha\vert \neq 1$ and $\vert a \vert =1$. If $\vert \alpha\vert = \vert a \vert =1$, 
then we see that the iterates $(g^n)_{n\in \Z}$ are contained in the compact group of transformations $(x,y)\mapsto (ux, vy)$ 
with $\vert u \vert = \vert v \vert =1$. So, $f$ is contained in a compact subgroup of $\Aut^0(X)$ and preserves a K\"ahler metric. 

\smallskip

{\bf{A3.--}} The last case is when $\vert \alpha \vert = 1$ and $\alpha^m=a$ for some $m \in \{0, 1, \ldots, d\}$. 
If we set $p(x)=\sum_{k\neq m} q_k (a-\alpha^k)^{-1}x^k$, we conjugate $f$ to $g_1(x,y)=(\alpha x, \alpha^m y + q_m x^m)$. Then, 
if we carry out a second conjugacy, over the open set $V:=\C^*\times \bbP^1(\C)$, by the map $h'(x,y)=(x, (q_m x^m)^{-1}y)$ we see
that $f$ is conjugate to $g_2(x,y)=(\alpha x, y+1)$ on $V$, with $\vert \alpha\vert =1$. This proves that $\hpol(f;V)=1$.

\smallskip

{\bf{B.--}} Now, assume that $f_B(x)=x+1$ in appropriate coordinates. We start by a lemma, which concerns the 
surface $X$ and its automorphism given by $(x,y)\mapsto (x+1, a y)$ in the chart 
$U:=\C\times \bbP^1$.

\begin{lem}
Let $X$ be, as above, a Hirzebruch surface of index $d$. Let $a$ be a complex number of modulus $1$, and
let $f$ be the automorphism of  $X$ defined by $f(x,y)=(x+1,ay)$ in $U=\C\times \bbP^1$. Then $\hpol(f; U)=0$.
\end{lem}

\begin{proof}
If $d=0$, $X$ is $\bbP^1(\C)^2$, and $\hpol(f)=0$ because the polynomial entropies of the linear projective transformations $x\mapsto x+1$ 
and $y\mapsto a y$ are equal to $0$ on $\C$ and $\bbP^1(\C)$ respectively (but $\hpol(f;X)=1$). 
There is a subtlety in the case of Hirzebruch surfaces. 
Fix a point $x$ in $\C$, and consider the sequence of points $(x+m, y)$ of $X$; in the second chart, we obtain the sequence 
$((x+m)^{-1}, (x+m)^dy)$; so, if we start with two points $(x,y)$ and $(x, y')$ which are $\epsilon$-close in the euclidean 
metric of $\C^2$, then after $m$ iterations of $f$, their images may be $\epsilon$-distinguished by looking at the second 
coordinate in the second chart. 

To control the polynomial entropy, we start by an example that concerns the case $d=1$.
Fix $\epsilon>0$, and $R=R(\epsilon)\geq 2$ such that every point $y$ of $\C\subset \bbP^1(\C)$ with $\vert y \vert > R$ (resp. $< 1/R$) is 
at distance $< \epsilon/2$ from $\infty$ (resp. $0$). Every integer $m$ between $0$ and $n$ 
can be written $m=2^k+b$ for some unique pair of non-negative integers $0\leq b <2^{k}$ and $0\leq k \leq \log_2(n)$. For each 
$k\leq \log_2(n)$, pick $N$ points $y_i$ such that the points $2^{k} y_i$ are $(\epsilon/2)$-dense in the annulus $1/R\leq \vert z \vert R$.
For a fixed $k$ we need $N\leq 4\pi (R/\epsilon)^2$ points, so that altogether we need at most $C(\epsilon) \log_2(n)$ points $y_i$.
Set $x_0=0$. Then, the points $(x_0, y_i)$ satisfy the following property: given any point $(x,y)$ with $\vert x- x_0\vert < \epsilon$  there is 
a point $y_i$ such that the points $((x+m)^{-1},(x+m) y)$ and $((x_0+m)^{-1}, (x_0+m) y_i)$ are $(2\epsilon)$-close for every $0\leq m \leq n$
in the metric of $\C\times \bbP^1(\C)$, and therefore also in the metric of $X$. 
Thus, the $n$-orbit of $(x,y)$ by $f$ is $\epsilon$-close to the $n$-orbit of $(0,y_i)$ in the Bowen metric. 

A similar result holds if we start with any point $x_j$ in place of the origin $x_0=0$ or if we change $d=1$ into any $d\geq 2$. 

Now, pick any compact subset $K$ of the open set $U=\C\times \bbP^1$, and fix $\epsilon > 0$ and $n\geq 2$. The
 projection of $K$ on $\C$ is contained in a compact 
disk $ \vert x \vert \leq r$, for some $r>1$; let $\{ x_j\; ; \; 1\leq j \leq J(r,\epsilon) \}$ be a finite $\epsilon$-dense set of points in that disk, 
with respect to the euclidean metric of $\C$ (we can take $J(r,\epsilon)\leq 2 \pi (r/\epsilon)^2$).
For each $j$, fix a subset $y_{j,i}$ of at most $C(\epsilon)\log(n)$ points as above. Altogether, we have $J(\epsilon,r)C(\epsilon)\log(n)$ 
points $(x_j, y_{j,i})$ and every point $(x,y)$ of $K$ is at distance $\leq 2\epsilon$ form one of these points in the Bowen metric $d_n^f$.
So, $\hpol(f;U)=0$.
\end{proof}

Let us now come back to the automorphism $f(x,y)=(x+1, ay + q(x))$ of $X$.
Consider the difference equation $ap(x)-p(x+1)=q(x)$;
if  $a\neq 1$, one sees by recursion on $d$ (starting with the highest degree terms) that there is a unique solution 
$p\in \C[x]$, of degree $\deg(p)=\deg(q)\leq d$. Then, the automorphism $h(x,y)=(x,y+p(x))$ of $X$ conjugates $f$ to 
$g(x,y)=(x+1, ay)$. 

If $a=1$, the difference equation becomes $p(x)-p(x+1)=q(x)$. One can find $p$ of degree $d$ such that $p(x+1)-p(x)=q(x)-q_d x^d$, 
so that $f$ is conjugate by an automorphism of $X$ to the transformation $g(x,y)=(x+1,y+q_d x^d)$. 
If $q_d\neq 0$, the unique fixed point of $g$ is the point $(\infty, \infty)$, and the dynamics in the complement of this fixed
point is wandering: as in the proof of Proposition~\ref{prop:pgl3}, the technics of Hauseux and Le Roux give $\hpol(g)=1$. 
If $q_d=0$, then $g(x,y)=(x+1,ay)$ with $a=1$. Thus, from the previous lemma we deduce that, in all cases, $\hpol(f)=1$ or 
$\hpol(f;U)=0$ (these two possibilities are not exclusive).

\subsubsection{Non minimal rational surfaces} 

Now, assume that $X$ is not a minimal rational surface, and consider an element $f$ of $\Aut(X)^0$. Fix a 
birational morphism $\eta\colon X\to X_0$ onto one of the minimal rational models of $X$. Since $\Aut(X)^0$ 
is connected, the lemma of Blanchard implies that there is an automorphism $f_0$ of $X_0$ such that 
$f_0\circ \eta= \eta \circ f$. Since $X_0$ is a minimal rational surface, it is isomorphic to $\bbP^2(\C)$, 
$\bbP^1(\C)\times \bbP^1(\C)$, or one of the Hirzebruch surfaces studied in Section~\ref{par:minimal_rational_surfaces}.
So, going from $(X,f)$ to $(X_0,f_0)$, and applying Section~\ref{par:minimal_rational_surfaces} or Proposition~\ref{prop:pgl3}
or Example~\ref{eg:quadric} we get a proof of Theorem~\ref{thm:rest-entropy} for all rational surfaces. With Section~\ref{par:Irrational_surfaces}, this 
concludes the proof of Theorem~\ref{thm:rest-entropy}.

\section{Automorphisms with dense orbits}\label{sec:no_finite_orbits}

In this section, we classify automorphisms of surfaces satisfying one of the following properties 
\begin{itemize}
\item all orbits of $f$ are infinite, i.e. $f\colon X\to X$ has no periodic orbit;
\item all orbits of $f$ are Zariski dense;
\item all orbits of $f$ are dense for the euclidean topology.
\end{itemize}
Those properties are listed from the weakest to the strongest: 
euclidean density implies Zariski density, which in turn excludes the existence of periodic orbit. 
The last property is exactly the notion of minimality with respect to the euclidean topology.
Our goal is to describe precisely, in dimension $2$, how the reverse implications fail. 
The results and proofs of this section are extremely close to~\cite{Reichstein-Rogalski-Zhang} (a reference we
discovered after writing this paper), the only difference being that we focus on the euclidean topology. 

\subsection{The Lefschetz formula and the Albanese morphism} 

\begin{lem}\label{lem:virtually-unipotent}
If $f$ is an automorphism of a compact K\"ahler surface without periodic orbits, its
 action on the cohomology of $X$ is virtually unipotent. 
\end{lem}

\begin{proof} If there is an eigenvalue $\lambda\in \C$ of $f^*$ with $\vert \lambda\vert >1$, we know from Section~\ref{sec:surfaces}
that $f$ has periodic orbits (see \cite{Cantat:Milnor} for a simple proof based on Lefschetz fixed point formula). 
We deduce that all eigenvalues of $f^*$ have
modulus $\leq 1$. Since $f^*$ preserves the integral cohomology $H^*(X;\Z)$, its eigenvalues are algebraic integers, and Kronecker
lemma shows that they are roots of unity. We deduce that some positive iterate $(f^n)^*$ is unipotent.
\end{proof}

Now, choose $n>0$ such that $(f^n)^*$ is unipotent; since $f^n$ has no fixed point, the holomorphic Lefschetz formula 
gives $h^{2,0}(X)-h^{1,0}+1=0$ (see~\cite{Cantat:Milnor}):

\begin{lem}\label{lem:h10}
If there is an automorphism of $X$ without periodic orbit, then  $X$ satisfies 
$
h^{1,0}(X)=h^{2,0}(X)+1.
$ 
In particular, there are non-trivial holomorphic $1$-forms on the surface $X$, and $X$ is not a rational surface. 
\end{lem}

Since $h^{1,0}(X)$ is positive, the Albanese map determines a non-trivial morphism $\alpha_X\colon X\to A_X$, where 
$A_X=H^0(X,\Omega^1_X)/H_1(X;\Z)$ is the Albanese torus. This map is equivariant with respect  $f$ 
and the automorphism $f_{alb}\colon A_X\to A_X$ induced by $f$:  this means that $f_{alb}\circ \alpha_X=\alpha_X\circ f$.

\subsection{Invariant genus $1$ pencil} 

Suppose that $f^*$ is unipotent and not equal to the identity. Then $f$ is parabolic and preserves a unique genus $1$ 
fibration $\pi\colon X\to B$ onto some Riemann surface $B$: there is an automorphism 
$f_B$ of $B$ such that $f_B\circ \pi = \pi\circ f$ (see Section~\ref{sec:surfaces}). Moreover, either $f_B$
 is periodic, or the surface $X$ is a torus
(see the proof of Proposition~\ref{pro:auto_surfaces_basics}). In particular, {\sl{if all orbits of $f$ are Zariski dense, then $X$ must be a torus}}. 

Assume  that  $f_B$ is periodic, and replace $f$ by $f^n$ where $n$ is the order of $f$ on the base. 
Given $b\in B$, $f$ preserves the fiber $X_b=\pi^{-1}(b)$. If $X_b$ is not a smooth curve of genus $1$, 
then there is a periodic orbit of $f$ in $X_b$. Thus, if all orbits of $f$ are infinite, then all 
fibers of $\pi\colon X\to B$ are smooth curves of genus $1$ (some of them may a priori be multiple fibers). 

\begin{lem}
Let $X$ be a compact K\"ahler surface, and let $f$ be an automorphism of $X$ with no finite orbit. 
Assume that $f^*$ is unipotent and not equal to the identity. Then $kod(X)\geq 0$, $f$ preserves a 
unique genus $1$ fibration, this fibration has no singular fiber (but may a priori have multiple fibers).
Moreover, $X$ is a minimal surface. 
\end{lem}

\begin{proof}
We already proved that $f$ preserves a genus $1$ fibration with no singular fiber. If $X$ were a 
ruled surface, then this ruling would be given by the Albanese map, $f$ would preserve two distinct
fibration, and this would contradict the fact that $f$ is parabolic. So, $kod(X)\geq 0$. Let us prove
the last assertion.

Let $\pi\colon X\to X_0$ be the birational morphism onto the minimal model $X_0$ of $X$. The exceptional divisor of $\pi$ 
coincides with the vanishing locus of all holomorphic $2$-forms of $X$. This implies that $f$ preserves this divisor, and some
iterate $f^n$ fixes each of its irreducible components. Those components being rational  curves, $f^n$ has a fixed point on 
each of them, contradicting the absence of periodic orbits. Thus, $X$ coincides with $X_0$. 
\end{proof}

\subsection{Tori and bi-elliptic surfaces} 

We consider two special cases, namely 
\begin{enumerate}
\item the minimal model of $X$ is a torus;
\item the minimal model of $X$ is a bi-elliptic surface. 
\end{enumerate}

If $X$ is a bi-elliptic surface, it is the quotient of  an abelian surface $A=B\times C$, where $B$ and $C$ are two elliptic curves, 
by a finite group $G$ acting diagonally on $A$: the action on $B$ is by translation $x\mapsto x+\epsilon$, and the action on $C$ is of the form 
\begin{equation}
y\mapsto \omega y+\eta\, , 
\end{equation}
where $\omega$ is a root of $1$ of order $2$, $4$, $3$ or $6$. The automorphism $f$ of $X$
lifts to an automorphism ${\tilde f}$ of the universal cover $\C^2$; here $\C^2=\C\times \C$, with coordinates $(x,y)$, and the elliptic
curves $B$ and $C$ are the quotients of the $x$-axis and the $y$-axis by lattices $\Lambda_B$ and $\Lambda_C$. Write
${\tilde f}(x,y)=L(x,y)+(a,b)$ for some linear transformation $L\in \GL_2(\C)$. Since ${\tilde f}$ covers $f$, the linear map normalizes the
linear part $(x,y)\mapsto (x,\omega y)$ of $G$. Thus, $L$ is a diagonal matrix. But from Lemma~\ref{lem:virtually-unipotent}, it is also 
virtually unipotent.  We deduce that $L^n={\sf Id}$ for some $n >0$. Changing $f$ in $f^{kn}$ for some $k>0$,
 we may assume that $f$ is covered by a translation that commutes to the linear part of $G$. Thus, {\sl some positive iterate $f^m$ 
 of $f$ is covered by a translation of type $(x,y)\mapsto (x+a,y)$}. This proves the following lemma. 
 
 \begin{lem}\label{lem:bi-elliptic}
 Let $f$ be an automorphism of a complex projective surface $X$ with no finite orbit.  
 If the minimal model of $X$ is bi-elliptic, then $X$ coincides with its minimal model, and a positive iterate of $f$ is covered by 
 a translation $(x,y)\mapsto (x+a,y)$ of a product $B\times C$ of two elliptic curves. In particular, the Zariski closure of each orbit
 is a curve of genus $1$.
 \end{lem}
 
 Let us now study the case of tori. 
 
 \begin{eg}
 Let $f$ be a translation on a $2$-dimensional compact torus $X=\C^2/\Lambda$, and let $M$ be the closure
 of the orbit of the neutral element $(0,0)\in X$. Then $M$ is a real Lie subgroup of $X$; its connected component 
 of the identity is a real torus $M^0\subset X$. The orbit $\{f^n(x);n\in \Z\}$ of any point $x\in X$ is dense in  $x+M$.
 Thus, on every compact torus, there are examples of translations whose orbits are Zariski dense but not dense for the euclidean topology. 
 \end{eg}
  
 \begin{eg}[Furstenberg~\cite{Furstenberg:1961}]
 Consider a $2$-dimensional torus $X$ which is the product of two copies of the same elliptic curve $E$; write $E=\C/\Lambda$ 
 and $X=\C^2/(\Lambda\times\Lambda)$ for some lattice $\Lambda\subset \C$. Then, consider the automorphism 
 \begin{equation}
 f(x,y)=(x+a, y+x+b)
 \end{equation}
 for some pair of elements $(a,b)\in E\times E$. Assume that $a$ is totally irrational with respect to $\C$, i.e. $x\mapsto x+a$ 
 has dense orbits in $E$ for the euclidean topology. Then, all orbits of $f$ 
 are dense for the euclidean topology. 
 \end{eg}
 
 \begin{lem}\label{lem:abel-surfaces}
 There are examples of automorphisms of complex abelian surfaces such that 
 \begin{itemize}
 \item $f$ is elliptic, all orbits of $f$ are Zariski dense, but no orbit is dense for the euclidean topology;
 \item $f$ is elliptic and all orbits of $f$ are  dense for the euclidean topology;
 \item $f$ is parabolic, and all orbits of $f$ are dense for the euclidean topology. 
 \end{itemize} 
If $f$ is an automorphism of an abelian surface with no finite orbit, then the following are equivalent: (i)
one orbit of $f$ is dense for the Zariski topology (resp. for the euclidean topology); (ii) every orbit of $f$ is dense for the Zariski topology (resp. for the euclidean topology).
 \end{lem}
 
 \begin{proof}
 Thanks to the previous examples, we only have to prove the second assertion. So, assume that there is an orbit of $f$ that is dense for the 
 euclidean topology; we want to prove that all orbits are dense. If some positive iterate of $f$ is a translation, this is easy. If not, 
 $f$ is parabolic; up to isogeny, we may assume that $X=E\times E$, with 
 \begin{equation}\label{eq:para-torus}
 f(x,y)=(x+a,y+kx)
 \end{equation}
 for some $k\geq 1$. Then, $x\mapsto x+a$ has dense orbits. If there is one orbit which is not dense, then 
 there is a non-trivial minimal invariant subset $M$ in $E\times E$. But this is impossible by Furstenberg's results. 
 
 Now, assume that there is an orbit of $f$ that is dense for the Zariski topology, and that $f$ is parabolic. Again, 
 $X=E\times E$, and we can assume that $f$ is as in Equation~\eqref{eq:para-torus}, with $x\mapsto x+a$ 
 a translation of infinite order. If the orbit of $(x,y)$ is not Zariski dense, then its Zariski closure is an $f$-invariant
 curve $C\subset E\times E$, on which $f$ induces an automorphism of infinite order. Thus, $C$ is 
 a curve of genus $1$, embedded in an $f$-invariant way into $E\times E$. But then, the translates of $C$ 
 form an $f$-invariant pencil, and since $f$ is parabolic, this pencil must coincide with the unique $f$-invariant
 fibration $(x,y)\mapsto x$. We get a contradiction because $a$ is not a torsion point of $E$. 
 \end{proof}

\subsection{Ruled surfaces (first step)}\label{par:RSFS}

Let us now assume that the Kodaira dimension of $X$ is $-\infty$. The Albanese map provides a fibration $\alpha_X\colon X\to A_X$, 
where $A_X=\C/\Lambda$ is a curve of genus $1$. The automorphism $f$ induces an automorphism $f_{alb}$ of $A_X$. If $f_{alb}^n(x)=x$
for some $m>0$, the fiber $\alpha_X^{-1}(x)$ is a curve of genus $0$ (it may be singular), and $f^m$ must fix a point in this fiber. 
Thus, the absence of finite orbit for $f$ implies that all orbits of $f_{alb}$ are infinite, and $f_{alb}$ is a translation of $A_X$ with Zariski dense
orbits. As a consequence, $\alpha_X$ is a submersion and $X$ is a fiber bundle over $A_X$ with rational fibers. The action of $f^*$ on 
$H^{1,1}(X;\R)$ can not be parabolic, because in that case $f$ preserves a unique fibration and this fibration is by curves of genus $1$. Thus, 
some positive iterate $f^m$ of $f$ is an element of $\Aut(X)^0$: there is a holomorphic vector field $\theta$ on $X$ and $f^m$ is the flow of $\theta$
at time $t=1$. This flow must permute the fibers of $\alpha_X$ and it is transverse to the fibration, exactly as 
in Section~\ref{par:minimal_ruled_inifinite}.

\begin{lem}\label{lem:kod-negative}
Let $f$ be an automorphism of a complex projective surface without finite orbits. 
If ${\sf kod}(X)$  is negative the albanese map $\alpha_X\colon X\to A_X$ is a submersion onto an elliptic curve whose
fibers are rational curves. Some positive iterate $f^m$ of $f$ is the flow, at time $1$, of a vector field which is everywhere
transverse to the fibration. 
\end{lem}

\begin{cor}\label{cor:kod-negative}
Under the same assumption on $f$ and ${\sf kod}(X)$, at least one
orbit of $f^m$ is contained, and Zariski dense, in an elliptic curve.
\end{cor}

\begin{proof} Consider the monodromy of the foliation induced by this vector field: it gives a representation of $\pi_1(A_X;x)$ in $\Aut(\alpha_X^{-1}(x))$, 
i.e. in  the group $\PGL_2(\C)$ of automorphisms  of $\bbP^1$. Since $\pi(A_X;x)\simeq \Z^2$, the monodromy group has a fixed point $p$. 
The orbit of $p$ under the flow of $\theta$ is a section of the fibration and is invariant under the action of $f^m$. Thus at least one
orbit of $f^m$ is contained, and Zariski dense, in an elliptic curve.\end{proof}

\subsection{Zariski dense orbits} 

\begin{thm}\label{thm:dense_torus}
 Let $f$ be an automorphism of a compact K\"ahler surface $X$.
If all orbits of $f$ are Zariski dense, then $X$ is a torus. On every torus there are translations whose orbits
are  dense for the euclidean and Zariski topologies.  
\end{thm}

\begin{que}
Let $f$ be an automorphism of a complex projective manifold $X$ of dimension $3$ (or more) acting minimally 
on $X$. Is $X$ automatically a torus ? 
\end{que}

\begin{proof}
If ${\sf kod}(X)=-\infty$, Corollary~\ref{cor:kod-negative} shows that $f$ has an orbit which is contained in 
a finite union of curves of genus $1$, contradicting our hypothesis. Thus  ${\sf kod}(X)\geq 0$: 
some positive multiple $mK_X$ of the canonical bundle has non-trivial sections.
Fix such a multiple, and consider the action of $f$ on the space of sections $H^0(X;mK_X)$. The existence of an eigenvector
provides a section $\omega$ of $mK_X$ such that $f^*\omega=\xi \omega$ for some $\xi \in \C^*$. The vanishing 
locus of $\omega$ is either empty, or an $f$-invariant curve. Since all orbits of $f$ are Zariski dense, $\omega$ does not vanish, $mK_X$ is the trivial bundle, and ${\sf kod}(X)=0$. 
Now, since  ${\sf kod}(X)=0$ and $h^{1,0}(X)>0$ (see Lemma~\ref{lem:h10}), the minimal model of $X$ is a torus or a bi-elliptic surface:
we conclude with Lemma~\ref{lem:bi-elliptic}. \end{proof}

\subsection{Ruled surfaces (second step)}\label{par:RSFS2} Let us come back to the study of ruled surfaces $\alpha_X\colon X\to A_X$ 
with an automorphism $f$ whose orbits are all infinite. According to Section~\ref{par:RSFS}, we can assume that $f\in \Aut(X)^0$
and $f$ is the flow of a vector field $\theta$ that is transverse to the fibration $\alpha_X$.

To simplify the notation, denote by $G$ the group $\Aut(X)^0$. By the universal property of the Albanese morphism, there is 
an $\alpha_X$-equivariant action of $G$ on $A_X$; and this action factors through the Abanese torus, via a homomorphism $A_G\to A_X$. 
Since the flow
$\Phi_\theta^t$ provides a non-trivial flow on $A_X$, we deduce that $G$ acts transitively on $A_X$. If $\dim(G)\geq 2$ the dimension of the kernel 
$H$ of the homomorphism $G\to A_X$ 
is positive,   $H$ has positive dimensional orbits in the fibers of $\alpha_X$, and $X$ is almost homogeneous: the group $G$ has an 
open orbit (for the Zariski topology). We shall study this case below. 

If $\dim(G) =1$, we have seen in Section~\ref{par:minimal_ruled_inifinite} that $G$ is an elliptic curve, isogeneous to $A_X$, and 
that no orbit of $f$ is Zariski dense (each orbit is contained in a $G$-orbit, hence in some elliptic curve).

To complete our study, we now rely on Section~3 of~\cite{Potters:1969} (namely the constructions on pages 251--253); see also Section~\ref{par:minimal_ruled_inifinite} above. 
Since $X$ is a ruled, almost homogeneous surface, $X$ is a topologically trivial $\bbP^1_\C$-bundle over the elliptic curve $A_X$, and there are two types of such bundles:
\begin{itemize}
\item[(a)] $X$ is the quotient of $\C^*\times \bbP^1_\C$ by the automorphism $(z,[x:y])\mapsto (\lambda z , [\mu x:y])$ for some 
pair $(\lambda, \mu)$ of complex numbers with $\lambda \mu\neq 0$ and $\vert \lambda\vert <1$.

\item[(b)] $X$ is the quotient of $\C^*\times \bbP^1_\C$ by the automorphism $(z,[x:y])\mapsto (\lambda z , [ x+y:y])$ for some complex number $\lambda$ 
with  $\vert \lambda\vert <1$. 
\end{itemize}

{\noindent{$\bullet$}} Case (a) corresponds in fact to two subcases. If $\mu$ is equal to $1$, then $X$ is just the product $A_X\times \bbP^1_\C$, and then 
every automoorphism is of the form $f\colon (z, [x:y])\mapsto (u(z), v[x:y])$ where $u$ is an automorphism of $A_X$ and $v$ is an element 
of $\PGL_2(\C)$. Thus, one gets: {\sl{No orbit of $f$ is dense for the euclidean topology; the general orbit of $f$ is Zariski dense if and only if $u$ and $v$ are two automorphisms of infinite order}}. If $\mu$ is a root of unity, the same result holds.

 Then, assume that $\mu$ is not a root of unity. 
Using affine coordinates $x=[x:1]$ for $\bbP^1_\C$, one sees that every automorphism
of $X$ comes from an automorphism of $\C^*\times \C^*$ of type $(z, x)\mapsto (\alpha z, \beta x)$ or $(\alpha z^{-1}, \beta x)$ for some
pair of complex numbers $(\alpha, \beta)$ with $\alpha \beta \neq 0$. Changing $f$ in $f^2$ we assume that $f(z, x)=(\alpha z, \beta x)$. 
Then, {\sl{the general orbit of $f$ is Zariski dense}}; indeed, the action of $f$ on $A_X$ has infinite order because otherwise $f$ has 
a periodic orbit, if the Zariski closure of a
general orbit is not $X$, then it is a multi-section of the fibration $\alpha_X$, but there are only two such multi-sections. A point $(z_0,x_0)$ has a
dense orbit for the euclidean topology if and only if for every point $(z_1,x_1)$ and every $\epsilon> 0$,   there are integers $m$ and $n$
such that $(\alpha^n z_0, \beta^n x_0)$ is $\epsilon$-close to $(\lambda^m z_1,\mu^m x_1)$ in $\C^*\times \C^*$. Taking logarithms, 
this means that the vectors $(\ln(\alpha), \ln(\beta))$, $(\log(\lambda),\log(\mu))$, $(2i\pi, 0)$, and $(0,2i\pi)$ generate a dense subgroup 
of $\C\times \C$, which of course is impossible since the rank of this group is at most $4$. This argument shows that in case (a), {\sl{there is
no automorphism with a dense orbit  for the euclidean topology}}. 

\smallskip
{\noindent{$\bullet$}} In case (b), every automorphism of $X$ can be written $(z,x)\mapsto (\alpha z , x+\beta)$ and, again, {\sl{the general orbit of $f$ is 
Zariski dense, but no orbit is dense for the euclidean topology}}. 

\subsection{Conclusion} 

\begin{thm}  
Let $f$ be an automorphism of a compact K\"ahler surface $X$, all of whose orbits are infinite. Replacing $f$
by some positive iterate, there are only three possibilities:

\smallskip

{\noindent}{\bf 1.--}   ${\sf kod}(X)=1$, $X$ is a fibration over a curve $B$ with smooth fibers of genus $1$ (some of them 
can be multiple fibers), and every orbit of $f$ is Zariski dense in such a fiber. 

\smallskip

{\noindent}{\bf 2.--}  ${\sf kod}(X)=0$, and  $X$ is a torus or  a bi-elliptic surface. If $X$ is bi-elliptic, then the Zariski closure of 
every orbit is a curve of genus $1$.  If $X$ is a torus, an orbit is dense for the Zariski (resp. for the euclidean) topology
if and only if all orbits are dense for this topology.

\smallskip

{\noindent}{\bf 3.--}      ${\sf kod}(X)=-\infty$, then $X$ is a ruled surface over an elliptic curve, at least one orbit is contained, and dense, 
in an elliptic curve (a section of the ruling), but no orbit is dense for the euclidean topology.  If the general orbit 
of $f$ is not Zariski dense, then

 (i) $X$ is isomorphic to the quotient of $\C^*\times \bbP^1_\C$
by 
$$(z,[x:y])\mapsto (\lambda z , [\xi x:y])
$$ 
with $0<\vert \lambda\vert <1$ and $\xi$ a root of unity, 

(ii) some positive iterate of 
$f$ is of the form $(z,[x:y])\mapsto (\alpha z, [\beta x:y])$ with $\beta$ a root of unity,

(iii)  the general orbits
are dense along multi-sections of the fibration~$\alpha_X$. 
\end{thm}

\begin{proof} All we have to do, is put together the previous results of this section together with Enriques-Kodaira classification 
of surfaces. First, if the Kodaira dimension of a projective variety is maximal, its group of automorphisms is finite. Thus,   ${\sf kod}(X)\leq 1$. 
Assume that ${\sf kod}(X)=1$. The action of $f$ on the base of the Kodaira-Iitaka fibration is periodic (see~\cite{Ueno}). 
Since the orbits of $f$ are infinite, every fiber of this fibration is a smooth curve of genus $1$, and every orbit is dense 
in such a fiber. 
When ${\sf kod}(X)=0$, we know that $X$ must be a minimal surface, and that $h^{1,0}(X)\geq 1$. Thus, $X$ is a torus
or a bi-elliptic surface. Then, we refer to Lemma~\ref{lem:bi-elliptic} and Lemma~\ref{lem:abel-surfaces}.
When ${\sf kod}(X)=-\infty$, we know from Section~\ref{par:RSFS} that $X$ is a ruled surface  over an elliptic curve, and
the conclusion follows from Section~\ref{par:RSFS2}. \end{proof}


\medskip
\begin{center}
{\bf{Part III.-- Small entropy and degree growth}}
\end{center}

\section{Small entropy: general facts}\label{sec:small_entropy}

In this section we gather a few remarks and examples concerning homeomorphisms of compact spaces
with small polynomial entropy. This section illustrates behaviors that may happen in the real analytic setting but do not 
occur for automorphisms. For example, Theorem \ref{thm:Borichev} provides an example of a diffeomorphism of a $2$-torus that has slow derivative growth and is not equicontinuous. This behavior never occurs for the automorphisms of smooth complex projective varieties by Theorem \ref{thm:derivative_bound}.
\subsection{Recurrence properties}

\begin{lem}\label{lem:obvious-small-entropy}
Let $f$ be a homeomorphism of a compact metric space $X$. 
\begin{enumerate}
\item If  $\hpol(f)< 1$, then for every $\epsilon >0$, there exists $k>0$ such that for every $x\in X$ there is a time $j\leq k$ with $dist(x,f^j(x))\leq \epsilon$; equivalently, all points of $X$ are recurrent, none of them is wandering. 
\item If the polynomial entropy
of $f$ is $< 1/2$, then $\limsup dist(f^n(x),f^n(y))\geq dist(x,y) /2$ for every pair of points
$(x,y)$ in $X\times X$.
\end{enumerate}
\end{lem}

\begin{proof}
The first part, due to \cite[Prop. 2.1]{Artigue-Carrasco-Monteverde}, is obtained as follows.  Suppose  there is $\epsilon >0$, 
such that for all $k>0$ one can find a point $y$ with $dist(y,f^j(y))>\epsilon$ for all $1\leq j\leq k$. Set $x_j=f^j(y)$
for $0\leq j\leq k$. Then, for the dynamics of $f^{-1}$, the points $x_j$ are $(\epsilon, k)$-separated, because if
$j<j'$ the distance between $f^{-j}(x_j)=y$ and $f^{-j}(x_{j'})=x_{j'-j}$ is greater than~$\epsilon$. Thus, $\hpol(f)\geq 1$. 
The equivalence with the recurrence property follows from the compactness of $X$.

The second assertion follows from the first one, applied to  $f \times f$.
\end{proof}

\begin{rem}
A point $x \in X$ is \emph{uniformly recurrent} for a homeomorphism $f:X \rightarrow X$ of a compact metric space if for any $\varepsilon>0$ there exists $N(\varepsilon)$ such that for any $n \in \N$ among any succesive iterates $f^{n+k}(x), k=0, \ldots, N-1$, there exists at least one such that $d(x, f^{n+k}(x))<\varepsilon$. Lemma~\ref{lem:obvious-small-entropy} does not say that all points are uniformly recurrent.
\end{rem}

\subsection{Growth of derivatives}

\begin{lem}\label{lem:lip}
Let $f$ be a homeomorphism of a compact manifold $X$. Denote by $Lip_f(n)$ 
the maximum of the Lipschitz constants of $\Id$, $f$, $\ldots$, $f^{n-1}$ ($Lip_f(n)$
is infinite if $f$ is not lipschitz). Then, $X$ is covered by  $O_\epsilon(Lip_n(f)^{\dim(X)})$ 
balls of radius $\leq \epsilon$ for the iterated metric $dist_n$. If $Lip(n)\leq n^\alpha$ for some
$\alpha >0$, 
then $\hpol(f)\leq \alpha \dim(X)$; if $Lip(n)=o(n^\alpha)$ for all $\alpha >0$, then $\hpol(f)=0$.
\end{lem}

\begin{proof}
If $x$ and $y$ satisfy $dist(x,y)\leq \frac{1}{2}\epsilon Lip_f(n)^{-1}$, then the distance between $f^k(x)$ and
$f^k(y)$ is less than $\epsilon$ for every natural integer $k \leq n-1$. And one can cover $X$ by 
roughly $(2\epsilon^{-1}Lip_f(n))^{\dim(X)}$ balls of radius $\frac{1}{2} \epsilon Lip_f(n)^{-1}$. 
\end{proof}

\begin{thm}
Let $f$ be a diffeomorphism of class ${\mathcal{C}}^2$ of a closed manifold $M$. Assume that the growth of
the derivative of $f^n$ is exponential: there is $\eta >0$ such that 
$
\parallel Df^n\parallel \geq \exp(\eta n)
$
as $n$ goes to $+\infty$. Then, $\hpol(f)\geq 1/2$. 
\end{thm}

For example, a diffeomorphism of the sphere with a north-south dynamics has polynomial entropy equal to $1$. It would be good to 
replace the inequality $\hpol(f)\geq 1/2$ by $\hpol(f)\geq 1$ in this theorem.

\begin{proof}[Sketch of the Proof]
There exists an $f$-invariant ergodic probability measure $\mu$ on $M$ with a positive Lyapunov exponent (See~\cite[\S~7.2]{Cantat:Bourbaki} for instance). 
Pesin's theory implies that a $\mu$-generic point $x$ has a non-trivial unstable manifold. Let $x$ and $y$ be points of 
such an unstable manifold. Then the distance between $f^n(x)$ and $f^n(y)$ goes to $0$ as $n$ goes to $-\infty$; by Lemma~\ref{lem:obvious-small-entropy}, 
this shows that $\hpol(f)\geq 1/2$. 
\end{proof}

\subsection{Skew products}

The following theorem answers a question of Artigue, Carrasco-Olivera, and Monteverde (see Problem~1 in~\cite{Artigue-Carrasco-Monteverde}). 
Let $\T$ denote the circle $\R/\Z$, so that $\T^d=\R^d/\Z^d$ is the torus of dimension $d$. 

\begin{thm}\label{thm:Borichev}
There exists a real analytic and area preserving diffeomorphism $f$ of the torus $\T^2$ satisfying the following four properties
\begin{enumerate}
\item $f$ is minimal;
\item its iterates $f^n$, $n\in \Z$, do not form an equicontinuous family;
\item for every $\epsilon >0$, the norm of the derivative satisfies $\parallel Df^n\parallel =o(n^\epsilon)$;
\item the polynomial entropy of $f$ vanishes.
\end{enumerate}
\end{thm}

\begin{rem}
Such examples exist on all tori $\T^k$, $k\geq 2$, but a homeomorphism $f$ of the circle 
with $\hpol(f)=0$ is conjugate to a rotation (see \cite{Labrousse:Circle}). 
\end{rem}

\begin{rem}[see \cite{Artigue-Carrasco-Monteverde}]
Let $\sigma$ be the shift on $\Lambda^\Z$ for some finite alphabet $\Lambda$.
The polynomial entropy of a subshift $\sigma_K\colon K\to K$ is $\geq 1$ for every 
$\sigma$-invariant infinite compact subset $K\subset \Lambda^\Z$. More generally, every expansive
homeomorphism of an infinite compact metric space has polynomial entropy $\geq 1$.
\end{rem}

 \begin{rem}\label{rem:Borichev}
Fix a function $\varphi\colon \R_+\to \R_+$ such that $\varphi$ is increasing, $\varphi$ is unbounded, $\varphi$ does not vanish, and
$\varphi(x)=o(x)$ as $x$ goes to $+\infty$. By a result of Borichev (see \cite{Borichev}, and also \cite{Polterovich:slow}), 
there is an analytic diffeomorphism $f$ of   $\R^2/\Z^2$
that  preserves the Lebesgue measure and satisfies
\begin{equation}
Lip_f(n) \leq \varphi(n) \quad {\text{and}} \quad \limsup_{n\to +\infty}\frac{Lip_f(n)}{\varphi(n)}>0.
\end{equation}
One can construct such an $f$ as a skew product $f(x,y)=(x+\alpha,y+g(x))$ for some well chosen periodic function $g$
and angle $\alpha$. The proof of Theorem~\ref{thm:Borichev} follows a similar strategy (and is simpler). 
 \end{rem}

 \subsubsection{Skew product} 
 Let $\alpha\in \R/\Z$
 be an irrational number, and let
 $g\colon \T\to \R$ be a continuous function such that $\int_0^1g(x)dx=0$. Consider the homeomorphism
 $f\colon \T^2\times \T^2$ defined by 
 \begin{equation}\label{eq:skew_product_on_torus}
 f(x,y)=(x+\alpha,y+g(x)).
 \end{equation}
 The $n$-th iterate of $f$ is $f^n(x,y)=(x+n\alpha, y+\sum_{j=0}^{n-1} g(x+j\alpha))$; if $g$ is smooth, $f$ is 
 a diffeomorphism, and the differential of $f^n$ is 
 \begin{equation}
 Df^n_{(x,y)}=\left( \begin{array}{cc} 1 & 0 \\ 
  \sum_{j=0}^{j=n-1}g'(x+j\alpha) & 1 \end{array}\right).
 \end{equation}
 
  \subsubsection{Minimality and equicontinuity} 
  
\begin{pro}[Furstenberg]\label{pro:Furstenberg} The homeomorphism $f$ is not minimal if and only if it is conjugate
to $(x,y)\mapsto (x+\alpha, y)$ by a homeomorphism $(x,y)\mapsto (x,y+h(x))$ with $h\colon \R/Z\to \R$ 
that solves the equation $h(x+\alpha)-h(x)=g(x)$. 
\end{pro}

This follows from \cite{Furstenberg:1961}. Indeed, Furstenberg proves that a proper minimal subset of the torus 
is the graph of such a homeomorphism $h$. 
   
  \begin{pro}\label{pro:Borichev1}
  If  $(f^k)_{k\in \Z}$ is an equicontinuous family, then 
  \begin{enumerate}
  \item  $\vert \sum_{j=0}^{n-1} g(x+j\alpha)\vert \leq B$ for some $B>0$ and all $n\geq 0$;
  \item $g$ is a coboundary: there is a continuous function $h\colon \R/\Z\to \R$ such that $h(x+\alpha)-h(x)=g(x)$ 
  for all $x\in \R/\Z$;
  \item $f$ is conjugate to $(x,y)\mapsto (x+\alpha,y)$ by a homeomorphism $(x,y)\mapsto (x,y+h(x))$;
  \item $f$ is not minimal.
  \end{enumerate}
  \end{pro}
  
This is well known to specialists, but we sketch the proof for completeness.
The family $f^\Z$ is equicontinuous; thus for any $\epsilon >0$ one can find $\eta>0$ such that $dist(f(x,y),f(x',y')) \leq \epsilon$
 as soon as $dist((x,y),(x',y'))\leq \eta$, where $dist$ is the euclidean distance on $\T^2$. Taking $\epsilon$ small, and covering $\T\times \{ 0\}$ by 
 $\eta^{-1}$ segments of length $\leq \eta$, one sees that  the image of $I\times \{ 0\}$ by any iterate $f^m$ of $f$ is a curve of
 length at most $\epsilon/\eta$. Now, take $B> \epsilon/\eta$. Assume that there exists $n$ 
 and $x$ with $\vert \sum_{j=0}^{n-1} g(x+j\alpha)\vert > B$. Since the mean
of $g$ is $0$, the Birkhoff sum $\sum_{j=0}^{n-1} g(x+j\alpha)$ vanish, and one can find an interval $I=[a,b]\subset \T$
such that $\sum_{j=0}^{n-1} g(a+j\alpha)=0$, $\sum_{j=0}^{n-1} g(x+j\alpha)>0$ (or $<0$) on $]a,b]$ and $\sum_{j=0}^{n-1} g(b+j\alpha)=B$ (or $-B$).
This implies that the segment $I\times \{ 0\}$ is mapped to a curve of length $\geq B$ by $f^n$, contradicting the choice of $B$. This
proves Assertion (1). 

Because the sums $\sum_0^{n-1} g(x+j\alpha)$ are uniformly bounded, and $x\mapsto x+\alpha$ is a minimal homeomorphism of $\T$, 
the lemma of Gottschalk and Hedlund (see \cite{Katok-Hasselblatt} page 100) shows that there exists a continuous function $h\colon \T\to \R$
satisfying $h(x+\alpha)-h(x)=g(x)$. This proves the second assertion, and the other two follow from it.

\subsubsection{Estimate of the derivatives}
 
Now, we choose $\alpha$ and $g$ explicitly. We will write $g$ as a Fourier series 
\begin{equation}\label{eq:analytic-g}
 g(x)=\sum_{k\in \Z} a_k e^{2i\pi k x}.
\end{equation} 
Fix a real number $r>1$. If  $a_k \leq r^{-k}$ then $g$ is an analytic function on the circle $\T$. 
To solve the equation $h(x+\alpha)-h(x)=g(x)$, we also expand $h$ as a Fourier series $\sum_{k} b_k e^{2i\pi k x}$; 
then, the $b_k$ must verify $b_k=(e^{2i\pi k\alpha}-1)^{-1} a_k$ for all $k\neq 0$. Choose 
\begin{equation}
\alpha = \sum_{i\geq 1} 10^{-q_i}=0.100010000000001000...
\end{equation}
where $q_1=1$, and the gaps $q_{n+1}-q_n$ between two consecutive $1$s increase quickly; more precisely, we shall 
assume that 
\begin{equation}\label{eq:choic-of-qn}
q_{n+1}-q_n > (\log(r)/\log(10)) 10^{q_n}.
\end{equation}
 Then, $10^{q_n}\alpha \simeq 10^{-(q_{n+1}-q_n)} \mod 1$.
This done, we choose $a_k=0$ for all indices except the one of the form $k=10^{q_n}$, in which case we 
choose $a_k=10^{-(q_{n+1}-q_n)}$. From Equation~\eqref{eq:choic-of-qn} we get $\vert a_k\vert \leq r^{-k}$ for all $k$, 
so that $g$ is analytic; but the solutions of the cohomological equation satisfy $b_k=1$ for $k\neq 1$, and we
deduce that there is no $L^2$ solution $h$ to the cohomological equation. This proves the first assertions of the following proposition. 

\begin{pro}\label{pro:Borichev2}
There is a pair $(\alpha, g)$ such that $\alpha$ is a Liouville number, $g$ is an analytic function,  the
cohomological equation $h(x+\alpha)-h(x)=g(x)$ ($\forall x \in \T$) has no continuous solution (resp. no $L^2$ solution), 
and the diffeomorphism $f$ satisfies $\parallel Df^n\parallel_{\T^2} =o(n^\epsilon)$ for all $\epsilon >0$. 
\end{pro}

Now, we want to find such a pair $(\alpha, g)$
satisfying 
\begin{equation}
\sum_{j=0}^{n-1} g'(x+j\alpha)=o(n^\epsilon)
\end{equation} 
for every $\epsilon >0$. To study this property, we expand $g$ in a Fourier series as in Equation~\eqref{eq:analytic-g}.
Fixing $n$, we set
\begin{equation}
D_n:=\sum_{j=0}^{n-1} g'(x+j\alpha)= 2i\pi \sum_{k\in \Z} \sum_{j=0}^{n-1} \left(e^{2i\pi k\alpha}\right)^j ka_k e^{2i\pi kx}
\end{equation}
and observe that 
 \begin{equation}
\vert \sum_{j=0}^{n-1} \left(e^{2i\pi k\alpha}\right)^j \vert \leq n \quad {\text{and}} \quad \sum_{j=0}^{n-1} \left(e^{2i\pi k\alpha}\right)^j=\frac{e^{2i\pi k\alpha n}-1}{e^{2i\pi k\alpha}-1}
\end{equation}
for all $n\geq 1$. Once $\epsilon$ has been fixed, we set $\tau=\epsilon /4 $ and split the sum $D_n$ in two parts:
\begin{equation}
\vert D_n\vert \leq 2\pi \left\vert \sum_{\vert k\vert \leq n^\tau} \left( e^{2i\pi k\alpha n} -1\right) \frac{k a_k}{e^{2i\pi k\alpha}-1} \right\vert + 2\pi \left\vert \sum_{\vert k \vert \geq n^\tau} n k a_k\right\vert.
\end{equation}
Since $g$ is analytic, its derivative is also analytic, and  $ka_k \leq C R^{-k}$ for some constants $C, R >1$. We shall assume that 
\begin{equation}
\left\vert \frac{a_k}{e^{2i\pi k\alpha}-1}\right\vert \leq 1,
\end{equation}
an inequality which is satisfied in the above construction of the pair $(\alpha, g)$. 
Altogether we get 
\begin{eqnarray}
\vert D_n\vert & \leq & 2\pi \times (  2n^\tau\times 2  n^\tau) + 2\pi \times 2\times \sum_{k\geq n^\tau} n C R^{-k}\\
& \leq & 16 \pi \times n^{2\tau} + 4\pi \times n \times C\frac{R}{R-1}R^{-n\tau}\\
& \leq & C' n^{2\tau}
\end{eqnarray}
because $nR^{-n\tau}\leq n^{2\tau}$ for $n$ large enough. Thus, $D_n\leq C'n^{\epsilon/2}=o(n^\epsilon)$, as required. 

\subsubsection{Conclusion}
The proof of Theorem~\ref{thm:Borichev} is now a direct consequence of Propositions~\ref{pro:Furstenberg},~\ref{pro:Borichev1} and~\ref{pro:Borichev2}, 
and Lemma~\ref{lem:lip}. 

\begin{que} Consider the family of skew products defined in Equation~\eqref{eq:skew_product_on_torus}: 
what is the value of $\hpol(f)$, as a function of $\alpha$ and of $g$ ?
\end{que}

\section{Slow growth  automorphisms}\label{sec:slow_growth}

In this section, we study automorphisms such that $\parallel Df^n\parallel $ grows slowly, or such that $f^\Z$ 
form an equicontinuous family on some large open subset of $X$, i.e. $f$ has a large Fatou component. 

\subsection{Grows of the derivative} The following result exhibits a growth gap phenomenon for the norm 
of the derivative of the iterates of any automorphism; from Remark~\ref{rem:Borichev}, we know that such a 
gap does not exist for real analytic diffeomorphisms of surfaces.

\begin{thm}\label{thm:derivative_bound}
Let $f$ be an automorphism of a smooth complex projective variety $X$ of dimension $d$. If 
$
\parallel Df^n\parallel = o(n)
$
then $f$ is an isometry of $X$ for some K\"ahler metric, and in particular the sequence $\parallel Df^n\parallel $ 
is bounded.
\end{thm}

\begin{lem}\label{lem:slow-auto-coho}
Let $f$ be an automorphism of a compact K\"ahler manifold $X$. Assume that the linear transformation 
$f^*\colon H^{1,1}(X;\R)\to H^{1,1}(X;\R)$ satisfies $\parallel (f^*)^n\parallel=o(n^2)$. Then some positive iterate $f^m$ of 
$f$ is in $\Aut(X)^0$. 
\end{lem}

\begin{proof} Let us show that $(f^*)^n$ is bounded 
on $H^{1,1}(X;\R)$.  By contradiction, writing $f^*$ in Jordan normal form, the only possibility would be that $(f^*)^n$ grows like $n$: all 
Jordan blocks of $f^*$ on $H^{1,1}(X;\R)$ have size at most $2$, with at least one of size $2$. 
Pick a vector $v$ in the interior of the K\"ahler cone. Then $\frac{1}{n} (f^*)^n(v)$ converges to a non-zero vector $w$ of the 
boundary of the K\"ahler cone as $n$ goes to $+\infty$, and at the same time $\frac{1}{n }(f^*)^{-n}(v)$ converges to the opposite vector $-w$. This  
 is impossible because the K\"ahler cone is salient, with a relatively compact base. 

Thus, the restriction of $(f^*)^n$ to $H^{1,1}(X;\R)$ is bounded;  the first item of the  Proposition 1.3.9 in  \cite{LoBianco}  shows that
$(f^*)^n$ is also bounded on $H^{p,p}(X;\R)$ for all $0\leq p\leq \dim(X)$, and then the second item shows 
that it is bounded on  $H^*(X;\R)$. Since it preserves the lattice $H^*(X;\Z)$, 
$f^*$ is periodic: there is a positive integer $n$ such that $(f^*)^n=\Id$ on $H^*(X;\R)$. 
Thus, by the theorem of Lieberman~(\cite{Lieberman}, Proposition~2.2 and Theorem~3.12), 
some positive iterate of $f$ is contained in $\Aut(X)^0$. \end{proof}

In view of Lemma~\ref{lem:slow-auto-coho}, Theorem~\ref{thm:derivative_bound} is now a direct consequence of the following proposition. 

\begin{pro}
Let $X$ be a complex projective manifold and $f$ be an element of $\Aut(X)^0$. If  
$
\parallel Df^n\parallel = o(n)
$
then $f$ is  in a compact subgroup of $\Aut(X)^0$.
\end{pro}

We thank Junyi Xie for the second step of the following proof, which is much simpler than our
initial strategy.

\begin{proof} Let $\alpha\colon X\to A_X$ be the Albanese map of $X$. There is a homomorphism $\rho\colon \Aut(X)^0\to \Aut(A_X)^0$ 
such that $\alpha\circ g=\rho(g)\circ \alpha$ for every $g\in \Aut(X)^0$. 

{\sl{Step 1. --}} Assume, first, that $\rho(f)=\Id$. Fix a very ample line bundle $L$
on $X$, and denote by $\psi_L\colon X\to \bbP^N(\C)$ the embedding given by the space of 
global sections of $L$, with  $N+1=\dim(H^0(X;L))$. 
Since $\rho(f)=\Id$, $f$ acts trivially on $\Pic^0(X)$ and  $f^*L=L$. In particular, $f$ acts linearly on
$H^0(X;L)$: this defines a linear projective
transformation $F\colon \bbP^N(\C)\to \bbP^N(\C)$ such that $\psi_L\circ f= F\circ \psi_L$. The
manifold $\psi_L(X)$ is not contained in a hyperplane; as a consequence, we can find a projective
basis $(x_1, x_2, \ldots, x_{N+2})$ of $\bbP^N(\C)$ whose elements $x_i$ are  in $\psi_L(X)$. Since $\max_{j\leq n} \parallel Df^n\parallel = o(n)$, 
we deduce that the distances $dist(F^n(x_i),F^n(x_j))$ are bounded from below by $\nu(n)/n$ for some 
sequence $\nu(n)$ that goes to $+\infty$ with $n$. This implies that $F$ is contained in a compact subgroup of $\Aut(\bbP^N(\C))\simeq \PGL_{N+1}(\C)$. Thus, $F$ preserves
some Fubini-Study form $\kappa$ on $\bbP^N(\C)$, and $f$ preserves $\kappa_X:=\psi_L^*(\kappa_{\psi_L(X)})$; thus, $f$ is 
contained in the compact subgroup of $\Aut(X)^0$ preserving $\kappa_X$. If $\rho(f)$ is an element of finite order, the same argument applies. 

\smallskip

{\sl{Step 2. --}}  Let us now assume that $\rho(f)$ has infinite order. Consider the Zariski closure $G$ of $f^\Z$ in $\Aut(X)^0$; 
changing $f$ in a positive iterate, we suppose that $G$ is an irreducible algebraic subgroup of the algebraic group $\Aut(X)^0$. 
Let $r$ be the complex dimension of $G$. As a real Lie group, $G$ is 
isomorphic to $\R^p/\Z^p\times \R^q$ for some pair of integers $(p,q)$ with $p+q=2r$. We endow $G$ with the riemannian 
metric given by the quotient of an euclidean metric on $\R^p\times \R^q$ for which $\R^p\perp \R^q$;  if $g$ is an element of $G$, 
then the sequence $g^n$ is bounded for the distance given by this metric if and only if $g$ is contained in the compact subgroup 
$\R^p/\Z^p\times \{ 0\}$. Consider the following subsets of $G$:
\begin{itemize}
\item $K=\{ g\in G\; ; \; (g^n)_{n\in \Z} \; {\text{is bounded in}} \; G\}$; this   is the maximal compact subgroup of $G$;
\item $J=\{ g\in G\; ; \; \parallel Dg^n\parallel = o(n) \}$, where $\parallel Dg\parallel$ is the maximum of the norm of 
$Dg_x\colon T_xX \to T_{g(x)}X$ for $x$ in $X$.
\end{itemize}
It suffices to prove that $J=K$. While $K$ is a subgroup, it is not clear that $J$ is a subgroup of $G$, but at least
we have $K\cdot J\subset J$ (and in particular $K\subset J$) because $K$ is compact.
From the first step, we deduce that $J\cap Ker(\rho)= K\cap Ker(\rho)$. So, we only need to prove that $\rho_{\vert K}\colon K\to \rho(G)$ 
is onto.

Firstly, $\rho(G)$ is a closed subgroup of $\Aut(A_X)^0\simeq A_X$ (acting by translations on $A_X$), so it is a compact
complex torus, isomorphic to $\R^{2g}/\Z^{2g}$ for some $g\geq 1$. Secondly, if $\rho_{\vert K}$ is not surjective, then 
$\rho$ induces a surjective homomorphism ${\overline{\rho}}\colon G/K\simeq \R^q\to \rho(G)/\rho(K)$ with $\rho(G)/\rho(K)$
a real compact torus of dimension $\geq 1$; but then the kernel of ${\overline{\rho}}$ and of $\rho$ would contain infinitely 
many connected components, in contradiction with the fact that $Ker(\rho)$ is an algebraic subgroup of the algebraic group $G$. 
So $\rho(K) = \rho(G)$. 
\end{proof}

\begin{thm}
Let $f$ be an automorphism of a compact K\"ahler manifold $X$. Assume that 
$\parallel Df^{n_i}_x\parallel\geq c n_i^{\rho}$
for all points $x\in X$, some infinite sequence of integers $n_i$, and some constant $c>0$. If $2\rho +1 > h^{1,1}(X)$,
then $\htop(f)>0$ and $\parallel Df^n\parallel \geq C\lambda^n$  for some constants $C>0$ and $\lambda >1$ and for all $n\geq 1$.
\end{thm}

Note that the assumption concerns $\parallel Df^n_x\parallel$ for every $x$, while the 
conclusion is on the supremum $\parallel Df^n\parallel$.
In the real analytic category, Herman constructs a diffeomorphism $f$ of a compact manifold $M$, of dimension $14$, such that
$\parallel Df^n_x\parallel\geq \exp(\sqrt{n})$ for all $n\geq 1$ and all $x\in M$ but $\htop(f)=0$ because
$\parallel Df^n\parallel$ grows sub-exponentially. See~\cite{Herman:NI}, ``{\sl{On a problem of Katok}}''.

\begin{proof} Fix a K\"ahler form $\kappa$ on $X$, and denote
by $d$ the dimension of $X$. Then $\int_X (f^n)^*\kappa\wedge \kappa^{d-1}\geq c_1 n^{2\rho}$ for some constant $c_1>0$. 
So, if $2\rho > h^{1,1}(X)-1$, we see that the linear transformation $f^*$ is not virtually unipotent on $H^{1,1}(X;\C)$: its spectrum contains 
an eigenvalue $\lambda\in \C$ of modulus $>1$; by Yomdin's theorem, the topological entropy of $f$ is $>0$; and $\parallel Df^n\parallel$
grows at least like $\sqrt{\vert \lambda\vert}^{n}$. 
\end{proof}

\subsection{A question}  

\begin{que}
Let $f$ be an automorphism of a compact K\"ahler manifold. If $\hpol(f)=0$ (resp. $\hpol(f)<1$), does it follow that 
$f$ is contained in a compact subgroup of $\Aut(X)$? 
\end{que}

We already gave a positive answer to this question when 
$\dim(X)\leq 2$ and when $X$ is a torus; the following proposition 
treats the case $f\in \Aut(X)^0$. 

\begin{pro} Let $X$ be a complex projective manifold.
If $f$ is an element of $\Aut(X)^0$ and $\hpol(f)<1$, then $f$ is in a 
compact subgroup of $\Aut(X)^0$. 
\end{pro}

\begin{proof} Denote by $G$ the Zariski closure of $f^\Z$ in $\Aut(X)^0$; changing $f$ in a 
positive iterate, we assume that $G$ is irreducible. For each $x\in X$, 
let $S_x\subset G$ be the stabilizer of $x$ in $G$. Then $G$ and $S_x$ are complex commutative
algebraic groups; the connected component $S_x^0$ has finite index in $S_x$. 
Below, we endow $G$ and $G/S_x$ with its (euclidean) topology of real Lie group.

\smallskip 

{\sl{Step 1. --}} 
Let $V(x)$ denote the Zariski closure of the orbit $G(x)\subset X$.
The algebraic group $G/S_x$ embeds as an open subset of $V(x)$. We say that a sequence $(g_n)$ 
of elements of $G$ goes to infinity in  $G/S_x$ if, given any compact subset $B$ 
of $G/S_x$ (for the euclidean topology), there is an integer $n_0$ such that $g_nS_X\notin B$ for $n\geq n_0$. 
If $(g_n)$ is such a sequence, then $g_n(x)$ goes to the boundary $V(x)\setminus G(x)$ of the orbit $G(x)$ 
when $n$ goes to infinity. Thus, if $(f^n)_{n\geq 0}$ goes to infinity in $G/S_x$, then $f^n(x)$ is a wandering 
orbit, and $\hpol(f)\geq 1$ by Lemma~\ref{lem:obvious-small-entropy}. 

Now remark that the following properties are equivalent in
$G/S_x$: $(f^n)_{n\geq 0}$ is unbounded, $(f^n)_{n\in \Z}$ is unbounded, and $(f^n)$ goes to infinity. 
We deduce that if $\hpol(f)=0$, then the class of $f$ modulo $S_x$ is contained in a compact subgroup of $G/S_x$ 
for every $x\in X$.

\smallskip 

{\sl{Step 2. --}}  Write $G=(\R^p/\Z^p)\times \R^q= \R^{p+q}/\Z^p$ as the quotient of the vector space $\R^{p}\times \R^q$ by 
the lattice $L=\Z^p$, and then $f=({\overline{u}}, v)$ as the projection of a vector $(u,v)\in \R^{p}\times \R^q$. 
A sequence $(g_n)$ goes to infinity in $G/S_x$ if and only if it goes to infinity in $G/S_x^0$, so we shall replace $S_x$ by $S_x^0$ in what follows. 
Then denote by $\Sigma_x^0$ the linear subspace of $\R^p\times \R^q$ whose projection in $G$ coincides with $S_x^0$.

From Step 1, $(f^n)$ is bounded in each of the quotients $A/S_x^0$:  for every $x\in X$, there is a point $w(x)\in \Sigma_x^0$
such that $(u,v)+w(x)\in \R^p$; this is equivalent to the inclusion $f\in (\R^p/\Z^p)$ modulo $S_x^0$. Now, consider an 
open subset $V$ of $X$ such that (1) $V$ is ${\mathcal{C}}^\infty$-diffeomorphic to $U\times \R^p/\Z^p$ for some 
open set $U\subset \R^k$ with $k+p=2\dim_\C(X)$ and (2) the action of $\R^p/\Z^p$ on $V$ is conjugate to the 
action by translation on the second factor of $U\times \R^p/\Z^p$. Such an open set exists, with $U$ a small transversal
of an orbit of $\R^p/\Z^p$ with minimal possible stabilizer. Our assumption on $f$ implies that $V$ is $f$-invariant 
and the action of $f$ on $V$ is conjugate to 
\begin{equation}
 (u,z) \in U \times \R^p/\Z^p \mapsto (u, z+\varphi(u)) \in U \times \R^p/\Z^p
\end{equation}
for some smooth map $\varphi\colon U\to   \R^p/\Z^p$. If $\varphi$ is not constant, then $\hpol(f)\geq 1$, and this 
contradicts our assumption. Thus, $\varphi$ is constant
and, on $V$, $f$ coincides with an element $g$ of $\R^p/\Z^p$. Since the action of $f$ and $g$ are (real) analytic
on $X$, we get $f=g$. Thus, $f$ is contained in a compact subgroup of $\Aut(X)$. \end{proof}

\subsection{Large Fatou components} 
The assumption in Theorem~\ref{thm:derivative_bound}  is global. In some situation, it is sufficient to study the iterates $f^n$
on some large open subset. 

Let $X$ be a compact projective manifold of dimension $d$, with an automorphism $f$. Embed $X$ in some projective space $\bbP^N(\C)$.
Denote by $H$ a hyperplane section of $X$, and by $\kappa$ the Fubini-Study (form restricted to $X$); the cohomology
class $[\kappa]$ is an element of the K\"ahler cone in $H^{1,1}(X;\R)$. 

\begin{lem} Let $U$ be an open set of $X$ such that 
\begin{itemize}
\item[(i)] $f^\Z$ form a normal family on $U$ (see Section~\ref{par:compact_kodaira});
\item[(ii)] $U$ contains a curve $C$ obtained by intersecting $d-1$ hyperplane sections: $C=H_1\cap \cdots \cap H_{d-1}$. 
\end{itemize}
Then, some positive iterate of $f$ is contained in $\Aut(X)^0$.
\end{lem}

\begin{proof} Since $C$ is a compact subset of $U$, the first hypothesis implies that the area of the curves $f^n(C)$ stays uniformly bounded. 
Equivalently, the norm of the cohomology classes $(f^n)^*[\kappa^{d-1}]$ stays bounded. On the other hand, the function
$f\mapsto \parallel f^*_{\vert H^{1,1}(X;\R)}\parallel$ and $f\mapsto \int_X (f^*\kappa)\wedge \kappa^{d-1}$ are comparable, 
because the K\"ahler cone is salient (with a relatively compact base) and its interior is non empty. So, we deduce that 
$\parallel (f^n)^*_{\vert H^{1,1}(X;\R)}\parallel$ is bounded; as a consequence, $f^*$ is contained in a compact subgroup 
$K$ of $\GL(H^{1,1}(X;\R))$ and $f^*$ preserves the k\"ahler class $\int_K g^* [\kappa] d{\mu}_K(g)$, where ${\mu}_K$
is the Haar measure of $K$ (normalized by ${\mu}_K(K)=1$). The conclusion follows from Lieberman's theorems. \end{proof}

Consider the set $Reg(f)\subset X$ of {\bf{regular points}} of $X$: $x\in Reg(f)$ if and only if one of the following
equivalent properties is satisfied:
\begin{itemize}
\item[(a)] there is a neighborhood $U$ of $x$ on which the iterates $f^n$, $n\in \Z$, form a normal family;
\item[(b)] for every $\epsilon>0$ there is a neighborhood $U$ of $x$ such that $dist(f^n(x),f^n(y))\leq \epsilon$
for all $n$ and all $y\in U$.
\end{itemize}
So, $f^\Z$ is equicontinuous on $V$ if and only if $V$ is contained in $Reg(f)$. Let $Irr(f)\subset X$ be 
the complement of $Reg(f)$, i.e. the set of irregular points. Note that we use the 
notion of regular points to match the vocabulary of \cite{Homma-Kinoshita}; equivalently, $Reg(f)$ 
is the Fatou set of $f$ (using both forward and backward iterates to define the Fatou set).

\begin{cor}\label{cor:reg_irr}
Let $X$ be a complex projective manifold. Let $f$ be an automorphism of $X$ such that 
$Reg(f)$ is connected  and contains a curve $C$ and $Irr(f)$ does not contain any algebraic curve. 
Then $f$ is contained in a compact subgroup
of $\Aut(X)$. \end{cor} 


\begin{proof} We can split $Reg(f)$ into two subsets: $Reg^{rec}(f)$ is the subset
of points $x\in Reg(f)$ such that $f^{n_i}(x)$ converges towards $x$ along a subsequence $n_i$ that
goes to $\infty$; $Reg^{wan}(f)$ is the complement of $Reg^{rec}(f)$ in $Reg(f)$. 

Apply Theorem~1 and~2 of \cite{Homma-Kinoshita} to $X$ and $f$, to get the following: 
either $Reg^{rec}(f)$ or $Reg^{wan}(f)$ is empty; if $Reg^{rec}(f)$ is empty then the limit set
of $(f^n(x))_{n\in \Z}$ is contained in $Irr(f)$ for all points in $Reg(X)$. 

Now, if $Reg^{rec}(f)$ is empty, we see that the sequence of curves $f^n(C)$ has bounded area
and converges towards a subset of $Irr(f)$; by Bishop theorem and the second assumption, we get
a contradiction.

So, $Reg(f)=Reg^{rec}(f)$ and $f$ is in $\Aut(X)^0$. As in the proof of Theorem~\ref{thm:derivative_bound}, denote 
by $A$ the Zariski closure of $f^\Z$ in $\Aut(X)^0$; changing $f$ into an iterate, we may assume
$A$ to be connected. For $x$ in $X$, let $S_x$ be the stabilizer of $x$ in $A$. 
Denote by  $U$ be the open subset of $X$ on which the stabilizer is minimal. The $f$-orbit of 
every $x$ in $U\cap Reg(f)$ is recurrent, so $f$ is contained in a compact subgroup of $A/S_x$. 
Since this holds for a general point of $X$, we conclude as in the proof of Theorem~\ref{thm:derivative_bound} 
that $f$ is in a compact subgroup of $\Aut(X)^0$. 
\end{proof}

\begin{cor}
Let $X$ be a connected complex projective manifold with  an automorphism $f$ for which
$Irr(f)$ is finite and non empty. Then $X$ is $\bbP^1(\C)$. 
\end{cor} 

\begin{proof}
If $\dim(X)>1$  there is a  curve in $Reg(f)$ and Corollary~\ref{cor:reg_irr}
provides 
a contradiction. So, $\dim(X)=1$ and $X$ is $\bbP^1(\C)$ because $Irr(f)$ is not empty.
\end{proof}

\section{Appendix: Proof of Lemma~\ref{lem:skew-P1}}

Changing $f$ into its inverse we may assume that $\vert a \vert > a_0$ for some real number $a_0>1$. 
Since every orbit is wandering, the polynomial entropy of $f$ is at least $1$ (see Example~\ref{eg:hpol1_homographies}). In the following we prove $\hpol(f)\leq 1$.

\smallskip

{\sl{Step 1.-- Preliminary remarks}}.
and $g(t)>t$ for every $t\notin \{0,1\}$, then there is an increasing homeomorphism $\varphi$ of $I$ that conjugates 
$h$ to $g$. Indeed, one can fix any increasing homeomorphism $\varphi_0$ from $[1/2, h(1/2)]$ to $[1/2, g(1/2)]$
and define $\varphi$ to be equal to $g^m\circ\varphi_0\circ h^{-m}$ on each of the intervals $h^m([1/2, h(1/2)]), m \in \Z$.

Moreover, if $H$ and $G$ are two homeomorphisms of $I^2$ of type 
\begin{equation}
H(s,t)=(s,h_s(t)), \quad G(s,t)=(s,g_s(t))
\end{equation}
where the homeomorphisms $h_s$ and  $g_s$ satisfy $h_s(t)>t$ and $g_s(t)>t$ for all $s \in I$ and all $t\notin \{0,1\}$, then the 
previous construction applied for every $s$ provides a homeomorphism $\Phi(s,t)=(s,\varphi_s(t))$ of $I^2$ that conjugates $H$ to $G$. 

\smallskip

Now, consider the euclidian space with coordinates $(u,v,w) \in \R^3$, and let $S\subset \R^3$ be the sphere  defined by the equation 
\begin{equation}
u^2+ v^2+ (w-1/2)^2= 1/4.
\end{equation}
Its center is $(0,0,1/2)$ and its radius is $1/2$. The points of $S$ can be parametrized by their height $w\in [0,1]$ and 
\textbf{longitudinal angle} $\theta\in [0,1]$, with $(u,v)=(w-w^2)^{1/2}(\cos(2\pi \theta), \sin(2\pi \theta))$.
By stereographic projection, the Riemann sphere $\bbP^1(\C)$ is homeomorphic to $S$, with $[0:1]$ corresponding
to the south pole $(0,0,0)$  and $[1:0]$ corresponding to the north pole $(0,0,1)$. 

Let $x=u+iv$ be a complex number of modulus $\rho$ and argument $2\pi \alpha(x)$. The preimage of the circle of radius $\rho$ centered at $(0,0)$ in the plane $\{w=0\}$ (containing $x$) under stereographic projection is the horizontal circle of $S$ at  height $w={\rho^2}/({\rho^2+1})$. Then the homography $[y_0:y_1]\mapsto [x y_0: y_1]$ is conjugate to the M\"obius transformation of $S$ given by $\theta\mapsto \theta+\alpha(x)$
and 
\begin{equation}
w\mapsto h_\rho(w)=\frac{\rho^2w}{(\rho^2-1)w+1}.
\end{equation}
If $\rho>1$, then $h_\rho(w)>w$ for all $w \notin \{0,1\}$.

\smallskip

{\sl{Step 2.--  A conjugacy}}. We first assume that $a$ is a holomorphic diffeomorphism onto its image. Then, one can use $a$ as a change of coordinate:
doing so, we replace $\overline\disk$ by a compact disk $\overline\Delta\subset \C$ contained in the open set $\{\vert x \vert >1\}$. 
We now assume that $a(x)=x$ for every $x\in \Delta$. 
Using the first step, we conjugate the map $f: \Delta \times\bbP^1(\C) \rightarrow \Delta \times \bbP^1(\C)$, 
\begin{equation}
f(x, [y_0:y_1])=(x, [x y_0: y_1]),
\end{equation}
to the homeomorphism $G_S: \Delta\times S \rightarrow \Delta\times S $
given by $G_S(x,w,\theta):=(x, g(w), \theta+\alpha(x))$; here $g$ is any homeomorphism of $[0,1]$ such that $g(w)>w$ if $w \notin\{0,1\}$, 
for instance one can take $g=h_\rho$ with $\rho=2$. In polar coordinates $x=(\rho, \alpha)$, this map $G_s$ 
becomes 
\begin{equation}
G_S(\rho, \alpha, w, \theta)= (\rho, \alpha, g(w), \theta+\alpha).
\end{equation} 
Since $x$ belongs to the compact 
disk ${\overline{\Delta}}\subset \{ \vert x\vert >1\}$,  $\rho$ and $\alpha$
stay in compact intervals $J_1$ and $J_2 \subset I$. Denote their lengths by $\ell_1$ and $\ell_2$ respectively. 

\smallskip

{\sl{Step 3.-- Polynomial entropy of $G_S$}}. Fix $\epsilon >0$, and let $N=N(\varepsilon)$ be the smallest positive integer such that $g^N(\sqrt{\epsilon}/3)>1-\sqrt{\epsilon}/3$. 
Then, there exists $\eta=\eta(\varepsilon) \in \R_{>0}$ such that for all $t,t' \in I$ with $t-t' \in (0,\eta)$, we have
\begin{equation}
\max_{0 \leq k \leq N} g^k(t)-g^k(t')\leq \sqrt{\epsilon}/3.
\end{equation}

Cover the interval $[\sqrt{\epsilon}/3, 1- \sqrt{\epsilon}/3]$ by $M$ intervals of length $< \eta$. The set of preimages by $g^{-\ell}$ for all $\ell\leq n$ of these intervals gives a cover of the interval $[g^{-n}(\sqrt{\epsilon}/3), 1- \sqrt{\epsilon}/3]$ into 
$m_0$ intervals with $m_0 \leq M(n+1)$. Then, choose a point $w_k$ in each of these intervals, thus defining the set $\{w_k\}_{k=1}^{m_0}$.
Now choose a set $\{r_j\}_{j=1}^{m_1}$ with $r_j \in J_1$ that $\epsilon$-covers $J_1$, and a set $\{\theta_l\}_{l=1}^{m_3}$ with $\theta_l \in I$ that $(\epsilon/2)$-covers $I$. We may impose $m_1 \leq \ell_1/\varepsilon+1$ and $m_3\leq 2/\varepsilon+1$. 

Given any integer $m\geq 4N/\epsilon$, we now construct a set $A(m)=\{\alpha_j\}_{j=1}^{m_2}$ with $\alpha_j \in J_2$, of size
$m_2\leq (4\ell_2 N)/\varepsilon+1$, that satisfies the following 
property : for every point $(\alpha, \theta)$, there is a point $(\alpha_j, \theta_l)$ such that for any $p=0, \ldots, N$ holds
\begin{equation}\label{eq:above}
\vert \theta + (m+p) \alpha - (\theta_l +(m+p)\alpha_j)\vert\leq \epsilon.
\end{equation}

Indeed, the left hand side of Equation \eqref{eq:above} is bounded from above by
\begin{equation}
  \vert (\theta + m \alpha - \theta_l) \vert + \vert p (\alpha-\alpha_j)-m\alpha_j\vert.
\end{equation}
There exists $l$ such that $\theta_l \in I$ is $\epsilon/2$-close to $\theta + m \alpha \in I$. Now we need to construct
a point $\alpha_j \in J_2$ such that 
\begin{equation}\label{eq:for_alpha}
\vert p (\alpha-\alpha_j) - m\alpha_j\vert\leq \epsilon/2.
\end{equation}
For this, we start with a set $\{\alpha_j'\}_{j=1}^{m_2}\subset  J_2$ that $(\epsilon/4N)$-covers $J_2$.  
Then, we perturb this choice into $\alpha_j=\alpha'_j+\nu_j$ with $\nu_j=\frac{\{m \alpha_j'\}}{m}$; here $\nu_j<\frac{\varepsilon}{4N}$ since $m\geq 4N/\epsilon$.
Now the left part of the Equation \eqref{eq:for_alpha} satisfies
\begin{equation}
\vert p (\alpha-\alpha_j')-pv_j-([m\alpha_j']+\{m \alpha_j'\}+m \nu_j)\vert \leq  \vert p (\alpha-\alpha_j') \vert + |p \nu_j|<\frac{\varepsilon}{2}.
\end{equation}

Once this is done, define $m(w)$ to be the smallest integer for which $g^m(w)\geq \sqrt{\epsilon}/3$.
There is a finite set ${\mathcal{S}}_0$ such that every point $(\rho, \alpha, w, \theta)$ with $m(w)\leq 4N/\epsilon$
is $\epsilon$-close to an element of ${\mathcal{S}}_0$ in the Bowen metric $d_n^{G_S}$ for every $n\geq 1$.
Now, by using the point sets we defined above, we define a final subset ${\mathcal{S}}(n;\epsilon) \subset \Delta \times S$ 
as follows: 
\begin{equation}
{\mathcal{S}}(n;\epsilon)= {\mathcal{S}}_0 \cup \{(r_i, \alpha_j, w_k, \theta_l) | m(w_k)>4N/\epsilon, \alpha_j \in A(m(w_k)) \}.
\end{equation}
This set is finite and its size is linear in $n$ (since $m_1, m_2$ and $m_3$ do not depend on $n$ and $m_0$ is linear in $n$).
By construction, every point $(\rho, \alpha, w, \theta) \in \Delta\times S$ is at distance less than $\epsilon$ from $ {\mathcal{S}}(n;\epsilon)$
with respect to the Bowen metric $d_n^{G_S}$. Indeed, for any $w \in I$, when $g^k(w)<\sqrt{\epsilon}/3$ or $g^k(w)> 1- \sqrt{\epsilon}/3$, the distances on the sphere $S$ are bounded by $\epsilon$. Moreover, for any orbit, there are at most $N$ iterates $g^s$
such that $ \sqrt{\epsilon}/3 < g^s(w) < 1- \sqrt{\epsilon}/3$ and these iterates correspond exactly to the values $s=m(w)+p$
with $p\leq N$. This shows that $\hpol(G_S)=1$ and hence $\hpol(f)=1$.

\smallskip

{\sl{Step 4.-- Conclusion}}.  Steps 2 and 3 show that the polynomial entropy of $f$ is equal to $1$ when $a$ is 
a holomorphic diffeomorphism onto its image. If the derivative of $a$ does not vanish, one covers $\overline\disk$ 
by smaller disks to reduce the proof to this latter case. If the derivative of $a$ vanishes, say at the origin, one can 
write $a(x)=a(0)+ x^d$ for some $d$ after a local change of the complex coordinate $x$. Then, the only change in 
the previous argument is that $\alpha$ is everywhere multiplied by $d$. This concludes the proof.

\begin{rem}
If we allow $\vert a(x)\vert=1$ for some values of $x\in \overline\disk$ the proof fails at the first and third steps. 
Indeed, if $h_s=\Id_{[0,1]}=g_s$ for some $s$, then the conjugacy $\Phi$ given by the $\varphi_s$ is not (at least not always)  
a homeomorphism of $I^2$. Morever, in the definition of $G_S$, 
one needs to take a family of homeomorphisms $g_{\rho(x)}$ depending on $\rho(x)$ such that $g_{\rho(x)}=\Id_{[0,1]}$ for some $x$;
when $\rho(x)$ approaches $1$, the time $N(x)$ for which $g_{\rho(x)}^N(\sqrt{\epsilon}/3)> 1-\sqrt{\epsilon}/3$ goes to $+\infty$, and the construction of the $\alpha_j$ breaks down.
 \end{rem}


%
%

%
%

\bibliographystyle{plain}
 
\bibliography{references}

\end{document}